%% file: actpaper.tex
\documentclass[submission,copyright,creativecommons]{eptcs}
\providecommand{\event}{SOS 2007} 
\usepackage{breakurl}             
\usepackage{underscore}           
\usepackage{amsthm}
\usepackage{amssymb}
\usepackage{xcolor}
\usepackage{mathtools}
\usepackage{tikzit}
\usepackage{stmaryrd}
\usepackage{tikz}
\usetikzlibrary{calc}
\usetikzlibrary{shapes.geometric}
\usetikzlibrary{shapes.misc}
\usetikzlibrary{positioning}

\input{sample.tikzstyles}

\definecolor{zx_grey}{RGB}{211,211,211}

\newcommand{\bvdots}{ \tikz[baseline, every node/.style={inner sep=0}]{ \node at (0,0){.}; \node at (0,-6pt){.}; \node at (0,6pt){.}; } }

\newtheorem{definition}{Definition}
\newtheorem{lemma}{Lemma}
\newtheorem{theorem}{Theorem}

\title{Propification and the Scalable Comonad}
\author{Titouan Carette
	\institute{Université Paris-Saclay, Inria, CNRS, LMF, 91190, Gif-sur-Yvette, France}
	\email{titouan.carette@universite-paris-saclay.fr}}
\def\titlerunning{Propification and the Scalable Comonad}
\def\authorrunning{Titouan Carette}

\begin{document}
	\maketitle
	
	\begin{abstract}
		String diagrams can nicely express numerous computations in symmetric strict monoidal categories (SSMC). To be entirely exact, this is only true for props: the SSMCs whose monoid of objects are free. In this paper, we show a propification theorem asserting that any SSMC is monoidally equivalent to a coloured prop. As a consequence, all SSMCs are within reach of diagrammatical methods. We introduce a diagrammatical calculus of bureaucracy isomorphisms, allowing us to handle graphically non-free monoids of objects. We also connect this construction with the scalable notations previously introduced to tackle large-scale diagrammatic reasoning.
	\end{abstract}
	
	\section*{Introduction}
	
	\noindent\textbf{Motivations:} The last two decades have seen the development of numerous graphical languages based on string diagrams. Among them, we mention in particular two lines of research that are still very active, Graphical linear algebra \cite{DBLP:journals/iandc/BonchiSZ17} for control flow graphs, and ZX-calculus \cite{DBLP:conf/icalp/CoeckeD08} for quantum computing. From a categorical point of view, those diagrammatic languages correspond to props, a particular kind of symmetric monoidal category (SMC), which are strict and have a free monoid of objects.
	
	The importance of strictness for string diagrammatic notations have been emphasised for a long time. The point is that when we represent the tensor product as parallel composition of diagrams, we don't want to consider any bracketing of the objects that would prevent us from composing by plugging wires. For example let's consider two arrows $f: (A\otimes B)\otimes C \to I$ and $g: I\to A\otimes (B\otimes C)$. Diagrammatically we would like to write:
	
	\begin{center}
		\input{stricdiag.tikz}
	\end{center}
	
	The two objects $(A\otimes B)\otimes C$ and $A\otimes (B\otimes C)$ have the same diagrammatical denotation but are not generally equal. So $f$ and $g$ cannot be composed. Diagrammatic notations cannot be used directly in a general SMC. However, there is a canonical way to put an adapter between the two morphisms to compose them. But then, we are not longer directly representing the arrows. MacLane's strictification theorem \cite{mac2013categories} ensures that any SMC is monoidally equivalent to an SSMC. Thus one can consider an equivalent SSMC before using the string diagrammatic notations. Sadly, this is still not enough. Even with an SSMC, we still can't, in general, directly use string diagrams. The tensor product on an SSMC induces a monoid structure on its set of objects. If this monoid structure is not free, we can find three objects $A$, $B$ and $C$ such that $A\otimes B = C$. If we consider two arrows $f:A\otimes B \to I$ and $g: I\to C$ we can represent them as:
	
	\begin{center}
		\input{freediag.tikz}
	\end{center}
	
	Here we are in the opposite situation as before. Nothing prevents us from composing those two diagrams in the category, but the combinatorial string diagrams cannot be plugged into each other. There is a diagrammatical typing problem. The only way to tackle those typing issues is to use a set of elementary wires whose tensor products freely generate all the other objects. In other words, the correct setting for string diagrams is props, not SSMCs. In addition to a strictification theorem, we need a propification theorem, and this is what this article proposes.\\
	
	\noindent\textbf{Contributions:} The main contribution of this paper is to identify a adjunction between \textbf{Prop} and \textbf{SSMC} inducing the scalable comonad on \textbf{Prop}. The left adjoint functor, called propification, turns any SSMC into a monoidally equivalent prop, allowing for manipulation of arrows in this category as string diagrams \textit{via} bureaucracy isomorphisms. Finally, relying on the work of \cite{wilson2022string}, we remark that the adjunction can be extended from $\textbf{SSMC}$ to the category $\textbf{SMC}$ on which the adjunction induces the usual strictification monad. The well-known strictification construction appears to be more than strictification: it is propification!\\
	
	\noindent\textbf{Related works:} This paper is the continuation of previous works on the foundation of scalable notations for large-scale diagrammatical reasoning. Those notations have been first introduced explicitly in \cite{chancellor2016graphical}, before being formalized in \cite{DBLP:conf/mfcs/CaretteHP19}. The formalism has then been refined and applied to quantum computing in \cite{carette2021quantum,DBLP:phd/hal/Carette21}. The fact that strictification factors through props have been first noticed by Paul Blain Levy, Sergey Goncharov, and Lutz Schroder \cite{traced}. They favoured the denomination concategories over colored props to emphasise the combinatorial nature of their definition and not to think of it as a special case of monoidal categories. Similar conclusions were reached independently by Paul Wilson, Dan Ghica and Fabio Zanasi \cite{wilson2022string} who recast the strictification process in a diagrammatical way. However, even if it is implicit in their work, the fact that strictification doesn't produce any SSMC but a prop is not explicitly stated.\\
	
	\noindent\textbf{Acknowledgment:} I thank Marc De Visme, who remarked that dividers and gatherers behave as the component of a natural transformation, and Paul Blain Levy, who kindly provided me more information on concategories. 
	
	\section{Background on props and Symmetric monoidal categories}
	
	This first section introduces the different definitions and notations that I will use throughout the paper.
	
	\begin{definition}
		A \textbf{small symmetric monoidal category} (SMC) $\mathbf{C}$ is a small category together with the following data:
		
		\begin{itemize}
			\item An object $I:\mathbf{C}$, the \textbf{tensor unit}.
			\item A \textbf{tensor product} functor $\_ \otimes \_ : \mathbf{C} \times \mathbf{C} \to \mathbf{C}$.
			\item Two natural isomorpisms with components $\lambda_c: I\otimes c \to c $, the \textbf{left unitor}, and $\rho_c : c\otimes I \to c$, the \textbf{right unitor}.
			\item A natural isomorphism with components $\alpha_{a,b,c}:  (a\otimes b)\otimes c  \to a\otimes (b\otimes c) $, the \textbf{associator}.
			\item A natural isomorphism with components $\sigma_{a,b}: a\otimes b \to b\otimes a$, the \textbf{swap}.
		\end{itemize}
		
		Furthermore, we require the following diagrams to commute:

		\begin{center}
			\input{symmoncomdiag.tikz}
		\end{center}
		
	\end{definition}
	
	A $SMC$ in which all $\alpha_{a,b,c}$, $\lambda_c $ and $\rho_c $ are identities is called a \textbf{symmetric strict monoidal category } (SSMC). There is a corresponding notion of functor:
	
	\begin{definition}
		A \textbf{strong symmetric monoidal} (SM) functor, is a functor $F:\mathbf{C}\to \mathbf{D}$ between two SMCs together with a natural isomorphim with components $\phi_{a,b}: F(a\otimes b)\to F(a)\otimes F(b) $ and an arrow $\phi: F(I) \to I$ such that the following diagrams commute:
		
		\begin{center}
			\input{symfunccomdiag.tikz}
		\end{center}
		
	\end{definition}
	
	A \textbf{strict symmetric monoidal} (SSM) functor is a SM functor for which all $\phi_{a,b}$ and  $\phi_c$ are identities. Finally there is also a corresponding notion of natural transformation
	
	\begin{definition}
		A natural transformation between two SM functors $\beta: F \Rightarrow G $ is \textbf{monoidal} if it satisfies:
		
		\begin{center}
			\input{monnatdiag.tikz}
		\end{center}
		
	\end{definition}
	
	Then two SM functors $F: \mathbf{C}\to \mathbf{D}$ and $G: \mathbf{D}\to \mathbf{C}$ form a \textbf{monoidal equivalence} if they form an equivalence of category where the natural isomorphisms $FG\Rightarrow id_{\mathbf{D}}$ and $GF\Rightarrow id_{\mathbf{C}}$ are also monoidal. We now have enough to state the key result behind strictifiaction.
	
	\begin{theorem}[MacLane Coherence Theorem]
		Given two functors $F,G:\mathbf{C}^n \to \mathbf{C} $ obtained from compositions of $I: 1 \to \mathbf{C}$ and $\_\otimes \_: \mathbf{C} \times \mathbf{C} \to \mathbf{C}$, there is a unique natural isomorphism $F\Rightarrow G$ that is made from compositions of $\alpha$, $\lambda$ and $\rho$.
	\end{theorem}
	
	A demonstration of the theorem is given in \cite{mac2013categories}. This result is known to be subtle to state precisely and to understand in detail. It is also very often distorted. Thus I advise the interested reader to consult the enlightening \cite{coh} for a more precise statement than the one given here. We now move to props.
	
	\begin{definition}
		A \textbf{prop} $\mathbf{P}$ is a SSMC such that the monoid structure induced by $\otimes$ on its set of object is free over a set of \textbf{colours} objects: $C(\mathbf{P})$.
	\end{definition}
	
	A \textbf{prop morphism} is then just an SSM functor between props. Notice that to define the action of a prop morphism on objects, we only need to define it on colours, and this is what we will do from now on.
	
	An object $\mathbf{a}:\mathbf{P}$ is always a finite list of colors $\sum\limits_{i} a_i $ where the concatenation of lists is denoted by $+$. The tensor unit is then the empty list denoted $0$. A colour is a one-element list. We will denote $C^* $ the set of finite lists over a set of colours $C$. We can picture an arrow $f: \sum\limits_{i=1}^n a_i \to \sum\limits_{j=1}^m b_j $ as a diagram with $n$ input wires and $m$ output wires. We should label those wires by the corresponding colour to be fully rigorous. We depict identities as simple wires and swaps as crossings; the identity of the empty list is the empty diagram.
	
	\begin{center}
		\input{stringdiagdef.tikz}
	\end{center}
	
	Those arrows can be composed as follow:
	
	\begin{center}
		\scalebox{0.8}{\input{stringdiagcomp.tikz}}
	\end{center}
	
	One can check that axioms of props correspond to intuitive, if not tautological, diagrammatic equations.
	
	\begin{center}
		\scalebox{0.8}{\input{stringdiageq.tikz}}
	\end{center}
	
	See \cite{selinger2010survey} for a review of the diverse results ensuring the soundness of those diagrammatical representations. Any prop can be defined by a set of generators and some equations between the compositions of those generators. More precisely, it has been shown in \cite{baez2017props} that the category $\mathbf{C}-\mathbf{Prop}$ of $C$-coloured props and prop morphisms is equivalent to the Eilenberg-Moore category of a monad over the functor category of signatures: $[C^* \times C^* ,\mathbf{set}]$. See \cite{DBLP:phd/hal/Carette21} for more details. Hence, a diagrammatical language is a prop axiomatised by generators and relations.
	
	\section{Propification}
	
	In this section we fix a symmetric strict monoidal category $\mathbf{C}$ and we describe how to construct from it a monoidally equivalent prop $\mathbf{Prop}(\mathbf{C})$. We denote $O(\mathbf{C})$ the set of objects of $\mathbf{C}$. A one element list will be denoted $[c]$ with $c:\mathbf{C}$.
	
	\begin{definition}
		The \textbf{content} $\{\mathbf{\ell}\}$ of a list $\mathbf{\ell}\in O(\mathbf{C})^* $ is defined inductively as: $\{0\} \doteq I $, $\{[c]\}=c$ and $ \left\{\mathbf{a}+\mathbf{b}\right\} \doteq \{\mathbf{a}\} \otimes \{\mathbf{b}\} $.
	\end{definition}
	
	This definition defines a content operator $\{\_\}: O(\mathbf{C})^* \to O(\mathbf{C}) $, as an example, we have $\left\{[a]+ [b] + [c]\right\}= a \otimes b\otimes c$.
	
	\begin{definition}[Propification]
		The prop $\mathbf{Prop}(\mathbf{C})$ is defined as the $O(\mathbf{C})$-colored prop such that:
		$\mathbf{Prop}(\mathbf{C})[\mathbf{a} ,\mathbf{b}]\doteq \mathbf{C}[\{\mathbf{a} \}, \{\mathbf{b} \}] $.
		Horizontal and sequential compositions, as well as swap maps are inherited from the symmetric monoidal structure of $\mathbf{C}$.
	\end{definition}
	
	There is a risk of confusion while considering arrows $f$ that can be seen either to be in $\mathbf{C}[\{\mathbf{a} \}, \{\mathbf{b} \}]$ or in $\mathbf{Prop}(\mathbf{C})[\mathbf{a} ,\mathbf{b}]$. To avoid this we will denote $[f]$ an arrow in $\mathbf{Prop}(\mathbf{C})[\mathbf{a} ,\mathbf{b}]$ corresponding to an arrow $f$ in $\mathbf{C}[\{\mathbf{a} \}, \{\mathbf{b} \}]$. Be careful, if $f$ denotes a unique arrow in $\mathbf{C}$ without ambiguity there are numerous arrows that can be denoted $[f]$ that differ by their types. They corresponds to the different input and output lists of objects having the same contents. Thus, I warn the reader that in the paper I will often use $[f]$ without more precisions. However, a quick type checking always allows to get rid of the ambiguity. The construction extends to a \textbf{propification functor}, but before defining it, we need to make some remarks on SM functors between SSMCs. Given a SM functor $F:\mathbf{C}\to \mathbf{D}$ we define for each list $\mathbf{c}\doteq\sum\limits_{i} [c_i ] \in O(\mathbf{C})^* $, an isomorphism $\phi_{\mathbf{c}}: F(\{\sum\limits_{i} [c_i ]\})\to \{\sum\limits_{i} [F(c_i ) ]\}$, defined by $\phi_{0}\doteq \phi$, $\phi_{[c]}\doteq id_{F(c)}$ and $\phi_{\mathbf{a}+\mathbf{b}}\doteq (\phi_{\mathbf{a}} \otimes \phi_{\mathbf{b}})  \circ \phi_{\{\mathbf{a}\},\{\mathbf{b}\}} $. The commutative diagrams in the definition of an SM functor ensure that $\phi_{\mathbf{c}}$ is uniquely defined. An important property is then that given $f:\{\mathbf{a}\}\to \{\mathbf{b}\}$ and $g:\{\mathbf{c}\}\to \{\mathbf{d}\}$ we have: $(\phi_{\mathbf{b}} \otimes \phi_{\mathbf{d}}) \circ (F(f)\otimes F(g)) \circ (\phi_{\mathbf{a}}^{-1}\otimes \phi_{\mathbf{c}}^{-1}) =\phi_{\mathbf{b}+\mathbf{d}} \circ F(f\otimes g) \circ \phi_{\mathbf{a}+\mathbf{c}}^{-1}$. We can now define our functor:
	
	\begin{align*}
		\mathbf{Prop}:\mathbf{SSMC}&\to \mathbf{Prop}& \mathbf{Prop}(F): \mathbf{Prop}(\mathbf{C})&\to \mathbf{Prop}(\mathbf{D})\\
		\mathbf{C} &\mapsto \mathbf{Prop}(\mathbf{C}) &  [c] &\mapsto [F(c)]\\
		F: \mathbf{C} \to \mathbf{D} & \mapsto \mathbf{Prop}(F): \mathbf{Prop}(\mathbf{C})\to \mathbf{Prop}(\mathbf{D}) & [f]:\mathbf{a}\to \mathbf{b}& \mapsto [ \phi_{\mathbf{a}} \circ F(f)\circ\phi_{\mathbf{b}}^{-1} ]
	\end{align*}
	
	Notice that here $\mathbf{Prop}(F)$ is defined only on colors $[c]$, from this we have that $\mathbf{Prop}(F)(\sum\limits_{i} [c_i ])=\sum\limits_{i} [F(c_i )] $. One can check that $\mathbf{Prop}(F)$ is well-typed and defines a prop morphism. There is an inclusion functor $\overline{\phantom{ }\cdot\phantom{ }}:\mathbf{Prop}\to \mathbf{SSMC}$, mapping a prop to the same prop seen as an SSMC and mapping prop morphisms to SSM functors. Given a prop $\mathbf{P}$, the corresponding SSMC category is denoted $\overline{\mathbf{P}}$. An important remark is that if we apply propification to a prop seen as an SSMC, we obtain a very different prop. In particular, the new colours are now lists of the former colours. We will look more closely at this phenomenon when we consider scalable construction. We define two monoidal functors:
	
	\begin{align*}
		\{\_\}_{\mathbf{C}}: \overline{\mathbf{Prop}(\mathbf{C})} &\to \mathbf{C} & [\_]_{\mathbf{C}} : \mathbf{C} &\to \overline{\mathbf{Prop}(\mathbf{C})} \\
		[c] &\mapsto c & c &\mapsto [c]\\
		[f]: \mathbf{a}\to \mathbf{b} &\mapsto f: \{\mathbf{a}\}\to \{\mathbf{b}\} & f: a\to b &\mapsto [f]: [a]\to [b]
	\end{align*}
	
	Note that $\{[\_ ]_{\mathbf{C}} \}_{\mathbf{C}} = id_{\mathbf{C}} $, and that $\{\_\}_{\mathbf{C}} $ is a strict monoidal functor while $[\_]_{\mathbf{C}} $ is only strong. In fact, $\{0\}_{\mathbf{C}} =I$ and $\{ \mathbf{a}+ \mathbf{b}\}=\{ \mathbf{a}\} \otimes \{\mathbf{b}\}$, while $[I]_{\mathbf{C}} \neq 0$ and $[a\otimes b]_{\mathbf{C}} \neq [a]_{\mathbf{C}} + [b]_{\mathbf{C}} $, however we have two natural isomorphisms $[id_I ]_{\mathbf{C}} :[I]_{\mathbf{C}} \to 0 $ and $[id_{a\otimes b}]_{\mathbf{C}} :[a\otimes b]_{\mathbf{C}} \to [a]_{\mathbf{C}} + [b]_{\mathbf{C}} $. Those two functors form a monoidal equivalence.
	
	\begin{theorem}[Monoidal equivalence] \label{thm:monequi}
		$\overline{\mathbf{Prop}(\mathbf{C})}$ is monoidally equivalent to $\mathbf{C}$. 
	\end{theorem}
	
	This theorem has important consequences for diagrammatic reasoning. Any SSMC can then be assigned a corresponding prop in which we can work diagrammatically. The following section describes how it works in practice.
	
	\section{Setting up an efficient diagrammatical bureaucracy}
	
	We now dive into the concrete diagrammatical formalism that allows drawing the arrows of an SSMC $\mathbf{C}$. The cornerstone is the introduction of bureaucratic isomorphisms that witness the equalities in $\mathbf{C}$ that no longer hold in $\mathbf{Prop}(\mathbf{C})$.
	
	\begin{definition}[Bureaucracy isomorphisms]
		An arrow $[f]:\mathbf{a}\to \mathbf{b}$ between lists $\mathbf{a}$ and $\mathbf{b}$ with similar content $c$, \textit{i.e.} such that $c=\{\mathbf{a}\}=\{\mathbf{b}\}$, is said \textbf{bureaucracy} if $f=id_{c}$. Note that a bureaucracy arrow is always an isomorphism.
	\end{definition}
	
	Those arrows are bureaucracy in the sense that they ensure the book-keeping of the equalities between tensor products in $\mathbf{C}$ that cannot holds in a prop. In a sense, we just relaxed those identities into isomorphisms.
	
	We have to make some brief remarks on bureaucracy isomorphisms. The identity arrows are bureaucracy isomorphisms, and bureaucracy isomorphisms are stable by inverse, tensor, and composition. Thus, we have a strict monoidal groupoid $\mathbf{Bur}(\mathbf{C})$ which is a subcategory of $\mathbf{Prop}(\mathbf{C})$ with the same objects. We can precisely describe $\mathbf{Bur}(\mathbf{C})$ as a disjoint union of cliques, one for each object of $\mathbf{C}$. Given an object $c:\mathbf{C}$, its clique as for vertices all formal tensors of objects in $\mathbf{C}$ that are equals to $c$ when evaluated in $\mathbf{C}$, in other words, all lists of objects in $\mathbf{C}$ with content $c$. The previous discussion can be summed up in the following Lemma:
	
	\begin{lemma}[Rewiring]
		Given two lists $\mathbf{a}, \mathbf{b}\in O(\mathbf{C})^* $ we have: $|\mathbf{Bur}(\mathbf{C})[\mathbf{a}, \mathbf{b}]|=\delta_{\{\mathbf{a}\}=\{\mathbf{b}\}}$.
	\end{lemma}
	
	The name \textit{rewiring} here borrowed from \cite{DBLP:conf/mfcs/CaretteHP19} will become clear with some string diagrams. We denote the unique arrow $\mathbf{a}\to \mathbf{b}$ in $\mathbf{Bur}(\mathbf{C})$ if it exists by the diagram:
	
	\begin{center}
		\input{bureaucracydef.tikz}
	\end{center}
	
	Then the rewiring lemma states that computing with bureaucracy isomorphisms is straightforward.
	
	\begin{center}
		\input{bureaucracycalc.tikz}
	\end{center}
	
	In fact, all bureaucracy morphisms can be obtain from the elementary dividers and gatherers, reminiscent of the ones in \cite{DBLP:conf/mfcs/CaretteHP19}, defined as:
	
	\begin{center}
		\input{divandgathdef.tikz}
	\end{center}
	
	and satisfying:
	
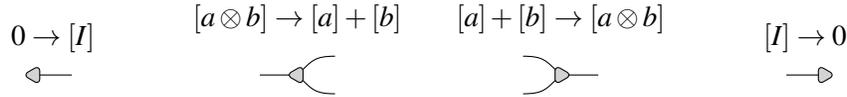
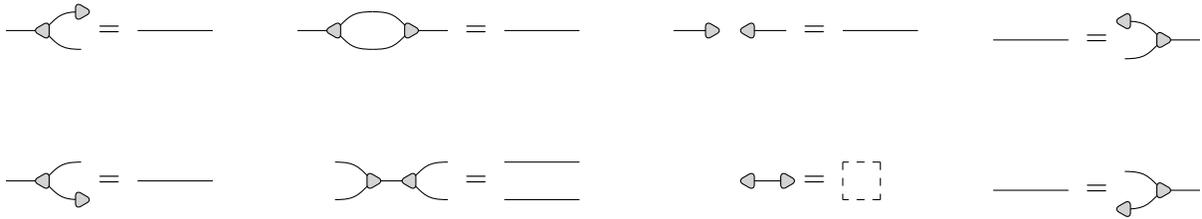
	\begin{center}
		\input{divandgatheq.tikz}
	\end{center}
	
	Graphically we represent the arrows $[f]$ as boxes indexed by $f$. The type is encoded in the inputs and output wires.
	
	\begin{center}
		\input{arrowbox.tikz}
	\end{center}
	
	Then, we can freely use the bureaucracy isomorphisms as adapters between the different incarnations in $\mathbf{Prop}(\mathbf{C})$ of the objects and arrows of $\mathbf{C}$.
	
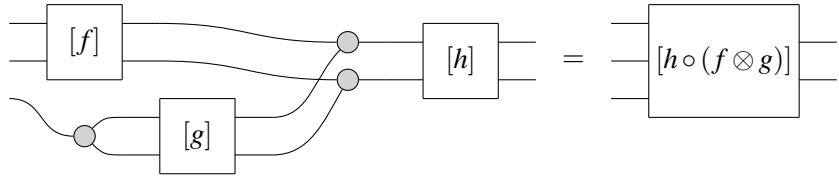
	\begin{center}
		\input{arrowboxex.tikz}
	\end{center}
	
	To end this section, let me present a concrete example of a category that could have been thought beyond the reach of diagrammatical technics but that we can tackle through propification. A typical prop used in diagrammatic quantum computing is $\mathbf{Qubits}$, obtained by taking as objects the integers and as arrows $n\to m$ the linear maps $\mathbb{C}^{2^n } \to \mathbb{C}^{2^m }$. Taking the usual tensor product of complex vector spaces, one can check that this is a prop in which many graphical languages have successfully been designed like the ZX-calculus \cite{DBLP:conf/icalp/CoeckeD08}. Let's consider a similar construction but replacing vector spaces with modules. We define $\mathbb{Z}\textbf{-mod}$ by again taking the integers as objects, but this time we take for arrows $n\to m$ the $\mathbb{Z}$-modules morphisms $\mathbb{Z}/n\mathbb{Z}\to \mathbb{Z}/m\mathbb{Z}$. With the conventions that $\mathbb{Z}/0\mathbb{Z} = \mathbb{Z}$ and $\mathbb{Z}/1\mathbb{Z} = \{0\}$. Considering the usual tensor product of $\mathbb{Z}$-modules, we have $\mathbb{Z}/n\mathbb{Z} \otimes_{\mathbb{Z}} \mathbb{Z}/m\mathbb{Z} \simeq \mathbb{Z}/\gcd(n,m)\mathbb{Z}$, and then in $\mathbb{Z}\textbf{-mod}$, $ n\otimes m = \gcd(n,m)$, the unit being $0$. The symmetry is taken to be the identity. Here, we have an example of a commutative strict monoidal category: a symmetric monoidal category whose swap is the identity. We see here that the monoid of objects is not free. Thus, drawing the arrows of such a category seems compromised. This is when propification comes into play. Then we can draw diagrams in $\mathbf{Prop}(\mathbb{Z}\textbf{-mod})$. Thus, we have now all the theoretical framework necessary to design graphical languages for $\mathbb{Z}$-modules.
	
	\section{The Scalable Comonad}
	
	The propification functor defined in the previous section has more uses than handling SSMCs pictorially. It is part of an adjunction inducing a comonad on $\mathbf{Prop}$ that corresponds to the scalable notations introduced in \cite{DBLP:conf/mfcs/CaretteHP19}.
	
	\begin{definition}
		The \textbf{scalable} functor is defined as: $S\doteq \mathbf{Prop}\circ \overline{\phantom{ }\cdot\phantom{ }}$, thus:
		
		\begin{center}
			\begin{align*}
				S:\mathbf{Prop}&\to \mathbf{Prop}& S(G): S\mathbf{P}&\to S\mathbf{Q}\\
				\mathbf{P} &\mapsto \mathbf{Prop}(\overline{\mathbf{P}}) &  [\sum\limits_{i} c_i ] &\mapsto [\sum\limits_{i} G(c_i) ]\\
				G: \mathbf{P} \to \mathbf{Q} &\mapsto S(G): S\mathbf{C}\to S\mathbf{D} & [f]& \mapsto [G(f)]
			\end{align*}
		\end{center}

	\end{definition}
	
	Given a prop $\mathbf{P}$, we define a prop morphism:
	
	\begin{align*}
		\langle\_\rangle_{\mathbf{P}}: S(\mathbf{P}) &\to \mathbf{P} \\
		[c] &\mapsto c \\
		[f]: \mathbf{a}\to \mathbf{b} &\mapsto f: \{\mathbf{a}\}\to \{\mathbf{b}\}
	\end{align*}
	
	Thus given a prop $\mathbf{P}$ and an SSMC $\mathbf{C}$, we now have two maps, $\langle\_\rangle_{\mathbf{P}}$ and $[\_]_{\mathbf{C}}$. They extend to natural transformations.
	
	\begin{lemma}\label{lm:naturality}
		There are two natural transformations: $[\_]: id_{SSMC} \Rightarrow \overline{\mathbf{Prop}(\_)}$ and $\langle\_\rangle: S \Rightarrow id_{\mathbf{Prop}}$ whose components are respectively: $[\_]_{\mathbf{C}}: \mathbf{C} \to \overline{\mathbf{Prop}(\mathbf{C})}$ and $\langle\_\rangle_{\mathbf{P}}:S(\mathbf{P})\to \mathbf{P}$.	
	\end{lemma}
	
	Those natural transformations form the unit and co-unit of an adjunction.
	
	\begin{lemma}\label{lm:adj}
		We have an adjunction:
		\begin{center}
			\input{ambadj.tikz}
		\end{center}
		
		In other words, given a prop $\mathbf{P}$ and a SSMC $\mathbf{C}$: $\overline{\left\langle [\_]_{\overline{\mathbf{P}}}\right\rangle_{\mathbf{P}}} = id_{\overline{\mathbf{P}}}$, and $\left\langle \mathbf{Prop}([\_]_{\mathbf{C}})\right\rangle_{\mathbf{Prop}(\mathbf{C})} = id_{\mathbf{Prop}(\mathbf{C})}$.	
	\end{lemma}

	This adjunction provides the endofunctor $S$ with a comonad structure. To explicit it, we will use the string diagramatic notations for natural transformations, see \cite{marsden2014category} for an introduction. The units and counits are depicted as follow:
	
	\begin{center}
		\begin{center}
			\input{unitcounitdef.tikz}
		\end{center}
	\end{center}
	
	We already know that they satisfies the equations:
	
	\begin{center}
		\input{unitcounitequations.tikz}
	\end{center}
	
	Then we can define a comonad structure on the functor $S$ as:
	
	\begin{center}
		\input{monadcomonaddef.tikz}
	\end{center}
	
	Note that we are not using those string diagrams in a rigorous way. The problem here is the same as with SSMCs. The types don't match graphically. The underlying (partial) monoid is not free. A rigorous approach would require the same kind of construction advocated in this paper but at the level of $2$-categories, which are generalisations of monoidal ones. However, we will not go this far here, even if I expect that very similar ideas would easily extend to this more general framework. Notice that structures very similar to bureaucracy isomorphisms are used in \cite{curien2008joy} to mimic equality in a $2$-categorical context. From the diagrams, we see that co-multiplication acts as:
	
	\begin{center}
		\begin{align*}
			\nu_{\mathbf{P}}: S\mathbf{P} &\to S^2 \mathbf{P} \\			
			\left[\sum\limits_{i} a_i\right]  &\mapsto \left[\sum\limits_{i} [a_i]\right] \\
			[f]&\mapsto [[f]]
		\end{align*}
	\end{center}
	
	Those structures are not completely new, they appeared numerous time in the development of scalable notations \cite{DBLP:conf/mfcs/CaretteHP19,carette2021quantum,DBLP:phd/hal/Carette21}. The idea is that while manipulating a large diagram, one sometimes would like to consider abstract wires representing an arbitrary number of small wires. This notation, which was informally used from the beginning of string diagram notations, is made formal by the scalable functor $S$. Starting with a prop $\mathbf{P}$ with wires coloured by a set $C$, we obtain a new prop $S(\mathbf{P})$  whose wires are coloured by lists of colours, conveniently representing groups of wires. The co-unit $\langle\_\rangle$ corresponds to the wire stripping functor of \cite{carette2021quantum}, it allows to open wires and split them into smaller ones. Such technics have so far mainly been applied to diagrammatical quantum computing \cite{DBLP:conf/mfcs/CaretteHP19,carette2021quantum}, but apply to any prop.
	
	\section{Beyond the strict case}
	
	Being mainly interested in direct applications to graphical language design and the connection to scalable notations, we restricted ourselves to the strict case of SSMCs. However, the construction we defined here generalises to non-strict SMCs. Then, instead of the lists of objects that are ubiquitous in this paper, we have to consider generalised lists with tensor units and bracketing. The content operator $\{\_\}$ is then defined to map those generalised lists to one specific bracketing that plays the role of a normal form. The propification must then be defined up to associators and unitors. However, it is still uniquely defined, thanks to MacLane's coherence Theorem.
	
	The reason why we do not carry those constructions in detail here, and only sketch the general idea, is that such a study has already been carried out extensively in \cite{wilson2022string}. There, the authors describe the strictification of SMCs in a diagrammatical way. Their construction is the straightforward generalisation of our propification functor to the non-strict case, and the result is also a prop. This is what we meant in the introduction by stating that the well-known strictification procedure dating from MacLane is, in fact, a propification procedure. So, we can substitute the usual adjunction between \textbf{SMC} and \textbf{SSMC}:
	
	\begin{center}
		\input{adjstrict.tikz}
	\end{center}
	
	with an adjunction between \textbf{SMC} and \textbf{Prop} that was already anticipated in \cite{traced}:
	
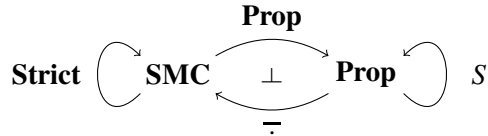
	\begin{center}
		\input{adjprop.tikz}
	\end{center}
	
	This adjunction induces the strictification monad on $\mathbf{SMC}$ as expected, but it is also interesting on the prop side, where it induces the scalable comonad. Finally, the $\textup{Prop}$ functor has an interest by itself as a way to apply string diagrammatic notations to any SMC. Then, in addition to bureaucracy isomorphisms, we have to consider new  isomorphisms playing the role of the associators and unitors in a similar way bureaucracy isomorphisms play the role of identities. For concrete examples, we redirect the reader to \cite{wilson2022string}. 
	
	\bibliographystyle{eptcs}
	\bibliography{generic}
	
	\appendix
	
	\section{Proofs}
	
	\begin{proof}[Proof of Theorem \ref{thm:monequi}]
		We define a natural transformation $\beta: id_{\overline{\mathbf{Prop}(\mathbf{C})}} \Rightarrow [\{\_\}_{\mathbf{C}} ]_{\mathbf{C}} $ with  components: $\beta_{\mathbf{c}}\doteq [id_{\{\mathbf{c}\}_{\mathbf{C}} } ]: \mathbf{c} \to  [\{\mathbf{c}\}_{\mathbf{C}} ]_{\mathbf{C}}$.
		
		We first check naturality, given an arrow $[f]:\mathbf{a}\to \mathbf{b}$ we have:
		
		\begin{center}
			$[\{[f]\}_{\mathbf{C}} ]_{\mathbf{C}} \circ \beta_{\mathbf{a}}= [f] \circ [id_{\{\mathbf{a}\}_{\mathbf{C}} }]=[f \circ id_{\{\mathbf{a}\}_{\mathbf{C}} }] =[id_{\{\mathbf{b}\}_{\mathbf{C}} }\circ f ]= [id_{\{\mathbf{b}\}_{\mathbf{C}} }]\circ [f]= \beta_{\mathbf{b}} \circ [f] $
		\end{center}
		
		The components of $\beta$ are indeed isomorphisms and we have:
		
		\begin{center}
			$\beta_{0}\circ id_{0}= [id_{\{0\}_{\mathbf{C}} } ]=[id_{I}]\quad$ and $\quad\beta_{\mathbf{a}+ \mathbf{b}}\circ id_{\mathbf{a}+ \mathbf{b}} = [id_{\{\mathbf{a}+ \mathbf{b}\}_{\mathbf{C}} } ]=[id_{\{\mathbf{a}\}_{\mathbf{C}} \otimes \{\mathbf{b}\}_{\mathbf{C}} }]$
		\end{center}

		$\beta$ is then a monoidal natural isomorphism between $id_{\overline{\mathbf{Prop}(\mathbf{C})}}$ and $[\{\_\}]$. Together with the fact that $\{[\_]_{\mathbf{C}} \}_{\mathbf{C}} = id_{\mathbf{C}} $ this gives us a monoidal equivalence between $\overline{\mathbf{Prop}(\mathbf{C})}$ and $\mathbf{C}$.
	\end{proof}
	
	\begin{proof}[Proof of Lemma \ref{lm:naturality}]
		
		We have to check naturality.
		
		\begin{center}
			\input{naturalityunitscounits.tikz}
		\end{center}
		
		Given $c:\mathbf{C}$, we have:
		
		\begin{center}
			$\overline{\mathbf{Prop}(F)}([c]_{\mathbf{C}}) = \overline{\mathbf{Prop}(F)}([c])= [F(c)]= [F(c)]_{\mathbf{D}}$
		\end{center}
		
		Given $f:a\to b$ an arrow in $\mathbf{C}$, we have:
		
		\begin{center}
			
			$\overline{\mathbf{Prop}(F)}([f]_{\mathbf{C}}) = \overline{\mathbf{Prop}(F)}([f])= [\phi_{[b]}\circ F(f)\circ \phi_{[a]}^{-1}]= [id_{[b]}\circ F(f)\circ id_{[a]}]= [F(f)]= [F(f)]_{\mathbf{D}}$
		\end{center}
		
		Given a color $\left[\sum\limits_{i} c_i \right]$ in $S(\mathbf{P})$, we have:
		
		\begin{center}
			
			$\left\langle S(G)\left(\left[\sum\limits_{i} c_i \right]\right)\right\rangle_{\mathbf{Q}} = \left\langle \left[G\left(\sum\limits_{i} c_i\right) \right]\right\rangle_{\mathbf{Q}}= G\left(\sum\limits_{i} c_i \right) = G\left(\left\langle\left[\sum\limits_{i} c_i \right]\right\rangle_{\mathbf{P}}\right) $
		\end{center}
		
		Given $[f]:\mathbf{a}\to \mathbf{b}$ an arrow in $S(\mathbf{P})$:
		
		\begin{center}
			$\left\langle S(G)\left(\left[f \right]\right)\right\rangle_{\mathbf{Q}} = \left\langle \left[G\left(f\right) \right]\right\rangle_{\mathbf{Q}}= G\left(f \right) = G\left(\left\langle\left[f \right]\right\rangle_{\mathbf{P}}\right) $
		\end{center}
		
	\end{proof}
	
	\begin{proof}[Proof of Lemma \ref{lm:adj}]
		Given an object $c:\mathbf{C}$, we have:
		
		\begin{center}
			$\left\langle \mathbf{Prop}([[c]]_{\mathbf{C}})\right\rangle_{\mathbf{Prop}(\mathbf{C})} = \left\langle [[c]]\right\rangle_{\mathbf{Prop}(\mathbf{C})} = [c]$
		\end{center}
		
		Given an arrow $[f]:\mathbf{a}\to \mathbf{b} $, we have:
		
		\begin{center}
			$\left\langle \mathbf{Prop}([[f]]_{\mathbf{C}})\right\rangle_{\mathbf{Prop}(\mathbf{C})} = \left\langle [\phi_{\mathbf{b}} \circ  [f]\circ \phi_{\mathbf{a}}^{-1}]\right\rangle_{\mathbf{Prop}(\mathbf{C})}= [id_{\{\mathbf{b}\}}] \circ  [f]\circ [id_{\{\mathbf{a}\}}]= [id_{\{\mathbf{b}\}} \circ  f\circ id_{\{\mathbf{a}\}}]=[f]$
		\end{center}
		
		Given an object $\sum\limits_{i} c_i : \mathbf{P}$, we have:
		
		\begin{center}
			$\overline{\left\langle \left[\sum\limits_{i} c_i \right]_{\overline{\mathbf{P}}}\right\rangle_{\mathbf{P}}} = \overline{\left\langle \left[\sum\limits_{i} c_i \right]\right\rangle_{\mathbf{P}}}= \sum\limits_{i} c_i  $
		\end{center}
		
		Given an arrow $f:\sum\limits_{i} a_i \to \sum\limits_{j} b_j $, we have:
		
		\begin{center}
			$\overline{\left\langle \left[f\right]_{\overline{\mathbf{P}}}\right\rangle_{\mathbf{P}}} = \overline{\left\langle \left[f \right]\right\rangle_{\mathbf{P}}}= f  $
		\end{center}
		
	\end{proof}
	
\end{document}

%% file: sample.tikzstyles
\tikzstyle{gn}=[fill=green, draw=black, shape=circle, tikzit category=ZX, tikzit fill=green, tikzit draw=black, tikzit shape=circle, inner sep=0.1em]
\tikzstyle{rn}=[fill=red, draw=black, shape=circle, tikzit fill=red, tikzit draw=black, tikzit category=ZX, tikzit shape=circle, inner sep=0.1em]
\tikzstyle{divide}=[regular polygon, regular polygon sides=3, shape border rotate=90, draw=black, fill={zx_grey}, inner sep=1.5pt, tikzit category=scal, rounded corners=0.8mm]
\tikzstyle{black}=[fill=black, draw=black, shape=circle, tikzit fill=black, tikzit draw=black, tikzit shape=circle, tikzit category=IH, inner sep=2pt]
\tikzstyle{gather}=[fill={zx_grey}, draw=black, tikzit category=scal, rounded corners=0.8mm, regular polygon, regular polygon sides=3, shape border rotate=-90, inner sep=1.5pt]
\tikzstyle{ggen}=[fill=white, draw=black, shape=rectangle, rounded corners=2mm, line width=1pt, tikzit draw=red, tikzit category=scal]
\tikzstyle{white}=[fill=white, draw=black, shape=circle, inner sep=2pt, tikzit category=IH]
\tikzstyle{mbox}=[fill=white, draw=black, rounded rectangle, rounded rectangle west arc=none, tikzit category=scal, tikzit shape=rectangle]
\tikzstyle{A}=[fill=white, shape=circle, tikzit category=scal, inner sep=1pt]
\tikzstyle{ggreen}=[fill=green, draw=black, shape=circle, tikzit category=SZX, tikzit fill=green, tikzit draw=black, line width=1pt, inner sep=0.1em]
\tikzstyle{gred}=[fill=red, draw=black, shape=circle, rounded corners=2mm, tikzit category=SZX, inner sep=0.1em, tikzit fill=red, line width=1pt]
\tikzstyle{ghad}=[minimum size=3mm, font={\scriptsize\boldmath}, shape=rectangle, inner sep=1mm, line width=1pt, outer sep=-1.5mm, scale=0.8, tikzit shape=rectangle, draw=black, fill=yellow, tikzit draw=blue]
\tikzstyle{boxm}=[fill=white, draw=black, rounded rectangle, tikzit category=scal, tikzit shape=rectangle, rounded rectangle east arc=none]
\tikzstyle{box}=[fill=white, draw=black, shape=rectangle]
\tikzstyle{had}=[fill=yellow, draw=black, shape=rectangle, tikzit category=ZX, tikzit fill=yellow, tikzit draw=black, inner sep=2.5pt]
\tikzstyle{gwhite}=[fill=white, draw=black, shape=circle, tikzit fill=white, tikzit shape=circle, line width=1 pt, inner sep=2 pt, tikzit draw=red]
\tikzstyle{gblack}=[fill=black, draw=black, shape=circle, tikzit fill=black, tikzit shape=circle, line width=1 pt, inner sep=2 pt, tikzit draw=red]
\tikzstyle{antipode}=[fill=red, draw=black, shape=rectangle, tikzit fill=red, tikzit draw=black, tikzit shape=rectangle, inner sep=2pt]
\tikzstyle{diamond}=[fill=white, draw=black, shape=diamond, inner sep=2pt]
\tikzstyle{mongr}=[fill=green, draw=green, shape=circle, inner sep=2pt]
\tikzstyle{monbl}=[fill=blue, draw=blue, shape=circle, inner sep=2pt]
\tikzstyle{bg}=[inner sep=0.7mm, minimum width=0pt, minimum height=0pt, fill=green, draw=white, very thick, shape=circle]
\tikzstyle{br}=[inner sep=0.7mm, minimum width=0pt, minimum height=0pt, fill=red, draw=white, very thick, shape=circle]
\tikzstyle{rmat}=[draw, signal, fill=red, signal to=east, signal from=west, inner sep=1pt, minimum height=6pt]
\tikzstyle{lmat}=[draw, signal, fill=red, signal to=west, signal from=east, inner sep=1pt, minimum height=6pt]
\tikzstyle{umat}=[draw, signal, fill=red, signal to=north, signal from=south, inner sep=1pt, minimum width=6pt]
\tikzstyle{dmat}=[draw, signal, fill=red, signal to=south, signal from=north, inner sep=1pt, minimum width=6pt]
\tikzstyle{box}=[shape=rectangle, text height=1.5ex, text depth=0.25ex, yshift=0.5mm, fill=white, draw=black, minimum height=5mm, yshift=-0.5mm, minimum width=5mm, font={\small}]
\tikzstyle{Z dot}=[inner sep=0mm, minimum size=2mm, shape=circle, draw=black, fill={rgb,255: red,160; green,255; blue,160}]
\tikzstyle{gdot}=[minimum size=3mm, font={\scriptsize\boldmath}, shape=rectangle, rounded corners=1.3mm, inner sep=1mm, outer sep=-1.8mm, scale=0.8, tikzit shape=circle, draw=black, fill=green, tikzit draw=blue]
\tikzstyle{X dot}=[Z dot, shape=circle, draw=black, fill={rgb,255: red,220; green,0; blue,0}]
\tikzstyle{rdot}=[minimum size=3mm, font={\scriptsize\boldmath}, shape=rectangle, rounded corners=1.3mm, inner sep=1mm, outer sep=-1.8mm, scale=0.8, tikzit shape=circle, draw=black, fill=red, tikzit draw=blue]
\tikzstyle{grdot}=[minimum size=3mm, font={\scriptsize\boldmath}, shape=rectangle, rounded corners=1.3mm, inner sep=1mm, line width=1pt, outer sep=-1.5mm, scale=0.8, tikzit shape=circle, draw=black, fill=red, tikzit draw=blue]
\tikzstyle{ggdot}=[minimum size=3mm, font={\scriptsize\boldmath}, shape=rectangle, line width=1pt, rounded corners=1.3mm, inner sep=1mm, outer sep=-1.5mm, scale=0.8, tikzit shape=circle, draw=black, fill=green, tikzit draw=blue]
\tikzstyle{arrow}=[-->]
\tikzstyle{rfarr}=[draw, signal, fill=black, signal to=east, signal from=west, inner sep=1pt, minimum height=6pt]
\tikzstyle{lfarr}=[draw, signal, fill=black, signal to=west, signal from=east, inner sep=1pt, minimum height=6pt]
\tikzstyle{ufarr}=[draw, signal, fill=black, signal to=north, signal from=south, inner sep=1pt, minimum width=6pt]
\tikzstyle{dfarr}=[draw, signal, fill=black, signal to=south, signal from=north, inner sep=1pt, minimum width=6pt]
\tikzstyle{ry}=[draw, signal, fill=yellow, signal to=east, signal from=west, inner sep=1pt, minimum height=6pt]
\tikzstyle{ly}=[draw, signal, fill=yellow, signal to=west, signal from=east, inner sep=1pt, minimum height=6pt]
\tikzstyle{uy}=[draw, signal, fill=yellow, signal to=north, signal from=south, inner sep=1pt, minimum width=6pt]
\tikzstyle{dy}=[draw, signal, fill=yellow, signal to=south, signal from=north, inner sep=1pt, minimum width=6pt]
\tikzstyle{new style 0}=[fill={zx_grey}, draw=black, shape=circle, inner sep=1mm]

\tikzstyle{arrow}=[->]
\tikzstyle{very thick}=[-, line width=1pt, tikzit draw=red]
\tikzstyle{reprise}=[-, line width=2pt, tikzit draw=blue]
\tikzstyle{pointille}=[dashed, -]
\tikzstyle{red}=[-, draw=red]
\tikzstyle{blue}=[-, draw=blue]
\tikzstyle{green}=[-, draw=green]
\tikzstyle{strike}=[-, tikzit draw={rgb,255: red,191; green,0; blue,64}, strike through]
\tikzstyle{strike'}=[-, tikzit draw=cyan, strike bend]
\tikzstyle{dashed arrow}=[->, tikzit draw=green, draw=black, dashed]

%% file: figures/stricdiag.tikz
\begin{tikzpicture}
	\begin{pgfonlayer}{nodelayer}
		\node [style=none] (0) at (0, 1.5) {};
		\node [style=none] (1) at (0, 0.5) {};
		\node [style=none] (2) at (3, 1.5) {};
		\node [style=none] (3) at (3, 0.5) {};
		\node [style=none] (4) at (0.5, 0.5) {};
		\node [style=none] (5) at (1.5, 0.5) {};
		\node [style=none] (6) at (2.5, 0.5) {};
		\node [style=none] (7) at (2.5, -0.5) {};
		\node [style=none] (8) at (1.5, -0.5) {};
		\node [style=none] (9) at (0.5, -0.5) {};
		\node [style=none] (10) at (1.5, 1) {$f$};
		\node [style=none] (11) at (14, -1.5) {};
		\node [style=none] (12) at (14, -0.5) {};
		\node [style=none] (13) at (11, -1.5) {};
		\node [style=none] (14) at (11, -0.5) {};
		\node [style=none] (15) at (13.5, -0.5) {};
		\node [style=none] (16) at (12.5, -0.5) {};
		\node [style=none] (17) at (11.5, -0.5) {};
		\node [style=none] (18) at (11.5, 0.5) {};
		\node [style=none] (19) at (12.5, 0.5) {};
		\node [style=none] (20) at (13.5, 0.5) {};
		\node [style=none] (21) at (12.5, -1) {$g$};
		\node [style=none] (22) at (1.5, -1.5) {$(A\otimes B) \otimes C $};
		\node [style=none] (23) at (12.5, 1.5) {$A\otimes (B \otimes C) $};
	\end{pgfonlayer}
	\begin{pgfonlayer}{edgelayer}
		\draw (0.center) to (1.center);
		\draw (1.center) to (3.center);
		\draw (3.center) to (2.center);
		\draw (2.center) to (0.center);
		\draw (4.center) to (9.center);
		\draw (5.center) to (8.center);
		\draw (6.center) to (7.center);
		\draw (11.center) to (12.center);
		\draw (12.center) to (14.center);
		\draw (14.center) to (13.center);
		\draw (13.center) to (11.center);
		\draw (15.center) to (20.center);
		\draw (16.center) to (19.center);
		\draw (17.center) to (18.center);
	\end{pgfonlayer}
\end{tikzpicture}

%% file: figures/freediag.tikz
\begin{tikzpicture}
	\begin{pgfonlayer}{nodelayer}
		\node [style=none] (0) at (0, 1.5) {};
		\node [style=none] (1) at (0, 0.5) {};
		\node [style=none] (2) at (2, 1.5) {};
		\node [style=none] (3) at (2, 0.5) {};
		\node [style=none] (4) at (0.5, 0.5) {};
		\node [style=none] (5) at (1.5, 0.5) {};
		\node [style=none] (8) at (1.5, -0.5) {};
		\node [style=none] (9) at (0.5, -0.5) {};
		\node [style=none] (10) at (1, 1) {$f$};
		\node [style=none] (11) at (11, -1.5) {};
		\node [style=none] (12) at (11, -0.5) {};
		\node [style=none] (13) at (10, -1.5) {};
		\node [style=none] (14) at (10, -0.5) {};
		\node [style=none] (17) at (10.5, -0.5) {};
		\node [style=none] (18) at (10.5, 0.5) {};
		\node [style=none] (21) at (10.5, -1) {$g$};
		\node [style=none] (22) at (1, -1.5) {$A\otimes B$};
		\node [style=none] (23) at (10.5, 1.5) {$C$};
	\end{pgfonlayer}
	\begin{pgfonlayer}{edgelayer}
		\draw (0.center) to (1.center);
		\draw (1.center) to (3.center);
		\draw (3.center) to (2.center);
		\draw (2.center) to (0.center);
		\draw (4.center) to (9.center);
		\draw (5.center) to (8.center);
		\draw (11.center) to (12.center);
		\draw (12.center) to (14.center);
		\draw (14.center) to (13.center);
		\draw (13.center) to (11.center);
		\draw (17.center) to (18.center);
	\end{pgfonlayer}
\end{tikzpicture}

%% file: figures/symmoncomdiag.tikz
\begin{tikzpicture}
	\begin{pgfonlayer}{nodelayer}
		\node [style=none] (0) at (4.5, -6.25) {$a\otimes b$};
		\node [style=none] (1) at (11.5, -6.25) {$a\otimes b$};
		\node [style=none] (3) at (5, -5.75) {};
		\node [style=none] (4) at (7.5, -3.25) {};
		\node [style=none] (5) at (8.5, -3.25) {};
		\node [style=none] (6) at (11, -5.75) {};
		\node [style=none] (7) at (5.5, -6.25) {};
		\node [style=none] (8) at (10.5, -6.25) {};
		\node [style=none] (9) at (5, -3.75) {$\sigma_{a,b}$};
		\node [style=none] (10) at (8, -2.75) {$b\otimes a$};
		\node [style=none] (11) at (8, -7) {$id_{a\otimes b}$};
		\node [style=none] (14) at (11, -3.75) {$\sigma_{b,a}$};
		\node [style=none] (15) at (-3.75, 6.25) {$a \otimes (I \otimes b)$};
		\node [style=none] (16) at (-7.25, 2.75) {$a\otimes b$};
		\node [style=none] (17) at (-8.75, 6.25) {};
		\node [style=none] (18) at (-10.75, 6.25) {$(a\otimes I)\otimes b $};
		\node [style=none] (19) at (-5.75, 6.25) {};
		\node [style=none] (20) at (-10.25, 5.75) {};
		\node [style=none] (21) at (-7.75, 3.25) {};
		\node [style=none] (22) at (-6.75, 3.25) {};
		\node [style=none] (23) at (-4.25, 5.75) {};
		\node [style=none] (24) at (-10.25, 3.75) {$\rho_a \otimes id_{b} $};
		\node [style=none] (25) at (-4.25, 3.75) {$id_{a} \otimes \lambda_{b} $};
		\node [style=none] (26) at (-7.25, 7) {$\alpha_{a,I,b} $};
		\node [style=none] (27) at (4, 5) {};
		\node [style=none] (28) at (8, 7.5) {$(a\otimes b) \otimes ( c\otimes d )$};
		\node [style=none] (29) at (3.5, 4.5) {$((a\otimes b) \otimes c)\otimes d $};
		\node [style=none] (30) at (12.5, 4.5) {$a\otimes (b \otimes ( c\otimes d ))$};
		\node [style=none] (31) at (7, 7) {};
		\node [style=none] (32) at (9, 7) {};
		\node [style=none] (33) at (12, 5) {};
		\node [style=none] (34) at (3.5, 3.5) {};
		\node [style=none] (35) at (3.5, 1.5) {};
		\node [style=none] (36) at (6.5, 0.5) {};
		\node [style=none] (37) at (9.5, 0.5) {};
		\node [style=none] (38) at (12.5, 1.5) {};
		\node [style=none] (39) at (12.5, 3.5) {};
		\node [style=none] (41) at (1.5, 2.5) {$\alpha_{a,b,c} \otimes id_d $};
		\node [style=none] (42) at (3.5, 0.5) {$(a\otimes (b \otimes c))\otimes d $};
		\node [style=none] (43) at (4.5, 6.5) {$\alpha_{a\otimes b,c,d} $};
		\node [style=none] (44) at (12.5, 0.5) {$a\otimes ((b \otimes c)\otimes d) $};
		\node [style=none] (45) at (14.5, 2.5) {$id_a \otimes \alpha_{b,c,d} $};
		\node [style=none] (46) at (11.5, 6.5) {$\alpha_{a, b,c\otimes d} $};
		\node [style=none] (47) at (8, -0.25) {$\alpha_{a, b\otimes c ,d} $};
		\node [style=none] (48) at (-10.75, -0.5) {$a \otimes (b\otimes c)$};
		\node [style=none] (49) at (-3.75, -0.5) {$(b\otimes c)  \otimes a$};
		\node [style=none] (50) at (-1.5, -3) {$b\otimes (c\otimes a)$};
		\node [style=none] (51) at (-13.25, -3) {$(a\otimes b )\otimes c$};
		\node [style=none] (52) at (-10.75, -5.5) {$(b\otimes a)\otimes c$};
		\node [style=none] (53) at (-3.75, -5.5) {$b\otimes (a\otimes c)$};
		\node [style=none] (54) at (-12.5, -2.5) {};
		\node [style=none] (55) at (-11.25, -1) {};
		\node [style=none] (56) at (-8.75, -0.5) {};
		\node [style=none] (57) at (-5.75, -0.5) {};
		\node [style=none] (58) at (-3.25, -1) {};
		\node [style=none] (59) at (-2.25, -2.5) {};
		\node [style=none] (60) at (-3.25, -5) {};
		\node [style=none] (61) at (-2.25, -3.5) {};
		\node [style=none] (62) at (-8.75, -5.5) {};
		\node [style=none] (63) at (-5.75, -5.5) {};
		\node [style=none] (64) at (-11.25, -5) {};
		\node [style=none] (65) at (-12.5, -3.5) {};
		\node [style=none] (66) at (-13.5, -1.5) {$\alpha_{a,b,c} $};
		\node [style=none] (67) at (-7.25, 0.25) {$\sigma_{a,b\otimes c} $};
		\node [style=none] (68) at (-13.5, -4.5) {$\sigma_{a,b}\otimes id_b  $};
		\node [style=none] (69) at (-7.25, -6.25) {$\alpha_{b,a,c} $};
		\node [style=none] (70) at (-1, -1.5) {$\alpha_{b,c,a} $};
		\node [style=none] (71) at (-1, -4.5) {$id_b \otimes \sigma_{a\otimes c} $};
	\end{pgfonlayer}
	\begin{pgfonlayer}{edgelayer}
		\draw [style=arrow] (7.center) to (8.center);
		\draw [style=arrow] (3.center) to (4.center);
		\draw [style=arrow] (5.center) to (6.center);
		\draw [style=arrow] (17.center) to (19.center);
		\draw [style=arrow] (20.center) to (21.center);
		\draw [style=arrow] (23.center) to (22.center);
		\draw [style=arrow] (27.center) to (31.center);
		\draw [style=arrow] (32.center) to (33.center);
		\draw [style=arrow] (34.center) to (35.center);
		\draw [style=arrow] (36.center) to (37.center);
		\draw [style=arrow] (38.center) to (39.center);
		\draw [style=arrow] (54.center) to (55.center);
		\draw [style=arrow] (56.center) to (57.center);
		\draw [style=arrow] (58.center) to (59.center);
		\draw [style=arrow] (60.center) to (61.center);
		\draw [style=arrow] (62.center) to (63.center);
		\draw [style=arrow] (65.center) to (64.center);
	\end{pgfonlayer}
\end{tikzpicture}

%% file: figures/symfunccomdiag.tikz
\begin{tikzpicture}
	\begin{pgfonlayer}{nodelayer}
		\node [style=none] (0) at (-3.5, 3) {$F(a\otimes b)$};
		\node [style=none] (1) at (3.5, 3) {$F(b\otimes a)$};
		\node [style=none] (2) at (-3.5, 0) {$F(a) \otimes F(b)$};
		\node [style=none] (3) at (3.5, 0) {$F(b)\otimes F(a)$};
		\node [style=none] (4) at (-3.5, 2) {};
		\node [style=none] (5) at (-3.5, 1) {};
		\node [style=none] (6) at (-1.5, 3) {};
		\node [style=none] (7) at (1.5, 3) {};
		\node [style=none] (8) at (-1.5, 0) {};
		\node [style=none] (9) at (1.5, 0) {};
		\node [style=none] (10) at (3.5, 2) {};
		\node [style=none] (11) at (3.5, 1) {};
		\node [style=none] (12) at (0, 3.75) {$F(\sigma_{a,b} )$};
		\node [style=none] (13) at (0, -0.75) {$\sigma_{F(a),F(b)}$};
		\node [style=none] (14) at (4.5, 1.5) {$\phi_{b,a}$};
		\node [style=none] (15) at (-4.5, 1.5) {$\phi_{a,b}$};
		\node [style=none] (16) at (12.5, -5) {$I \otimes F(c)$};
		\node [style=none] (17) at (19.5, -5) {$F(I)\otimes F(c)$};
		\node [style=none] (18) at (19.5, -8) {$F(I\otimes c)$};
		\node [style=none] (19) at (12.5, -8) {$F(c)$};
		\node [style=none] (20) at (12.5, -6) {};
		\node [style=none] (21) at (12.5, -7) {};
		\node [style=none] (22) at (14.5, -5) {};
		\node [style=none] (23) at (17.5, -5) {};
		\node [style=none] (24) at (17.5, -8) {};
		\node [style=none] (25) at (14.5, -8) {};
		\node [style=none] (26) at (19.5, -6) {};
		\node [style=none] (27) at (19.5, -7) {};
		\node [style=none] (28) at (16, -4.25) {$\phi \otimes id_{F(c)}$};
		\node [style=none] (29) at (16, -8.75) {$F(\lambda_{c})$};
		\node [style=none] (30) at (20.5, -6.5) {$\phi_{I,c}^{-1}$};
		\node [style=none] (31) at (11.5, -6.5) {$\lambda_{F(c)}$};
		\node [style=none] (32) at (12.5, 3) {$F(c) \otimes I$};
		\node [style=none] (33) at (19.5, 3) {$F(c)\otimes F(I)$};
		\node [style=none] (34) at (19.5, 0) {$F(c\otimes I)$};
		\node [style=none] (35) at (12.5, 0) {$F(c)$};
		\node [style=none] (36) at (12.5, 2) {};
		\node [style=none] (37) at (12.5, 1) {};
		\node [style=none] (38) at (14.5, 3) {};
		\node [style=none] (39) at (17.5, 3) {};
		\node [style=none] (40) at (17.5, 0) {};
		\node [style=none] (41) at (14.5, 0) {};
		\node [style=none] (42) at (19.5, 2) {};
		\node [style=none] (43) at (19.5, 1) {};
		\node [style=none] (44) at (16, 3.75) {$\phi \otimes id_{F(c)}$};
		\node [style=none] (45) at (16, -0.75) {$F(\rho_{c})$};
		\node [style=none] (46) at (20.5, 1.5) {$\phi_{c,I}^{-1}$};
		\node [style=none] (47) at (11.5, 1.5) {$\rho_{F(c)}$};
		\node [style=none] (48) at (-5, -4) {$F(a) \otimes (F(b)\otimes F(c))$};
		\node [style=none] (49) at (5, -4) {$F(a)\otimes F( b \otimes c)$};
		\node [style=none] (50) at (7.25, -6.5) {$F(a\otimes (b\otimes c))$};
		\node [style=none] (51) at (-7.5, -6.5) {$(F(a)\otimes F(b) )\otimes F(c)$};
		\node [style=none] (52) at (-5, -9) {$F(a\otimes b)\otimes F(c)$};
		\node [style=none] (53) at (5, -9) {$F((a\otimes b)\otimes c)$};
		\node [style=none] (54) at (-6.75, -6) {};
		\node [style=none] (55) at (-5.5, -4.5) {};
		\node [style=none] (56) at (-1.5, -4) {};
		\node [style=none] (57) at (1.5, -4) {};
		\node [style=none] (58) at (5.5, -4.5) {};
		\node [style=none] (59) at (6.5, -6) {};
		\node [style=none] (60) at (5.5, -8.5) {};
		\node [style=none] (61) at (6.5, -7) {};
		\node [style=none] (62) at (-1.5, -9) {};
		\node [style=none] (63) at (1.5, -9) {};
		\node [style=none] (64) at (-5.5, -8.5) {};
		\node [style=none] (65) at (-6.75, -7) {};
		\node [style=none] (66) at (-8.5, -5) {$\alpha_{F(a),F(b),F(c)} $};
		\node [style=none] (67) at (0, -3.25) {$id_{F(a)} \otimes \phi_{b,c}^{-1}$};
		\node [style=none] (68) at (-8.5, -8) {$\phi_{a,b}^{-1}\otimes id_{F(c)} $};
		\node [style=none] (69) at (0, -9.75) {$\phi_{a\otimes b,c}^{-1} $};
		\node [style=none] (70) at (8.5, -5) {$\phi_{a,b\otimes c}^{-1} $};
		\node [style=none] (71) at (8.5, -8) {$F(\alpha_{a,b,c})$};
	\end{pgfonlayer}
	\begin{pgfonlayer}{edgelayer}
		\draw [style=arrow] (6.center) to (7.center);
		\draw [style=arrow] (4.center) to (5.center);
		\draw [style=arrow] (10.center) to (11.center);
		\draw [style=arrow] (8.center) to (9.center);
		\draw [style=arrow] (22.center) to (23.center);
		\draw [style=arrow] (20.center) to (21.center);
		\draw [style=arrow] (26.center) to (27.center);
		\draw [style=arrow] (24.center) to (25.center);
		\draw [style=arrow] (38.center) to (39.center);
		\draw [style=arrow] (36.center) to (37.center);
		\draw [style=arrow] (42.center) to (43.center);
		\draw [style=arrow] (40.center) to (41.center);
		\draw [style=arrow] (54.center) to (55.center);
		\draw [style=arrow] (56.center) to (57.center);
		\draw [style=arrow] (58.center) to (59.center);
		\draw [style=arrow] (60.center) to (61.center);
		\draw [style=arrow] (62.center) to (63.center);
		\draw [style=arrow] (65.center) to (64.center);
	\end{pgfonlayer}
\end{tikzpicture}

%% file: figures/monnatdiag.tikz
\begin{tikzpicture}
	\begin{pgfonlayer}{nodelayer}
		\node [style=none] (0) at (0, 1.5) {$F(a\otimes b)$};
		\node [style=none] (1) at (7, 1.5) {$G(a\otimes b)$};
		\node [style=none] (2) at (0, -1.5) {$F(a) \otimes F(b)$};
		\node [style=none] (3) at (7, -1.5) {$G(a)\otimes G(b)$};
		\node [style=none] (4) at (0, 0.5) {};
		\node [style=none] (5) at (0, -0.5) {};
		\node [style=none] (6) at (2, 1.5) {};
		\node [style=none] (7) at (5, 1.5) {};
		\node [style=none] (8) at (2, -1.5) {};
		\node [style=none] (9) at (5, -1.5) {};
		\node [style=none] (10) at (7, 0.5) {};
		\node [style=none] (11) at (7, -0.5) {};
		\node [style=none] (12) at (3.5, 2.25) {$\beta_{a\otimes b}$};
		\node [style=none] (13) at (3.5, -2.25) {$\beta_a \otimes \beta_b $};
		\node [style=none] (14) at (8, 0) {$\phi_{a,b}$};
		\node [style=none] (15) at (-1, 0) {$\phi_{a,b}$};
		\node [style=none] (33) at (17, 1.5) {$I$};
		\node [style=none] (34) at (19.5, -1) {$G(I)$};
		\node [style=none] (35) at (14.5, -1) {$F(I)$};
		\node [style=none] (36) at (15, -0.5) {};
		\node [style=none] (37) at (16.5, 1) {};
		\node [style=none] (40) at (16, -1) {};
		\node [style=none] (41) at (18, -1) {};
		\node [style=none] (42) at (19, -0.5) {};
		\node [style=none] (43) at (17.5, 1) {};
		\node [style=none] (45) at (17, -1.75) {$\beta_I $};
		\node [style=none] (46) at (19, 0.5) {$\phi$};
		\node [style=none] (47) at (15, 0.5) {$\phi$};
	\end{pgfonlayer}
	\begin{pgfonlayer}{edgelayer}
		\draw [style=arrow] (6.center) to (7.center);
		\draw [style=arrow] (4.center) to (5.center);
		\draw [style=arrow] (10.center) to (11.center);
		\draw [style=arrow] (8.center) to (9.center);
		\draw [style=arrow] (36.center) to (37.center);
		\draw [style=arrow] (42.center) to (43.center);
		\draw [style=arrow] (40.center) to (41.center);
	\end{pgfonlayer}
\end{tikzpicture}

%% file: figures/stringdiagdef.tikz
\begin{tikzpicture}
	\begin{pgfonlayer}{nodelayer}
		\node [style=none] (0) at (1, 1) {};
		\node [style=none] (1) at (1, -1) {};
		\node [style=none] (2) at (3, 1) {};
		\node [style=none] (3) at (3, -1) {};
		\node [style=none] (4) at (0, 0.5) {};
		\node [style=none] (5) at (0, -0.5) {};
		\node [style=none] (6) at (3, 0.5) {};
		\node [style=none] (7) at (4, 0.5) {};
		\node [style=none] (8) at (1, -0.5) {};
		\node [style=none] (9) at (1, 0.5) {};
		\node [style=none] (10) at (2, 0) {$f$};
		\node [style=none] (17) at (4, -0.5) {};
		\node [style=none] (18) at (3, -0.5) {};
		\node [style=none] (24) at (0.5, 0) {$\bvdots$};
		\node [style=none] (25) at (3.5, 0) {$\bvdots$};
		\node [style=none] (26) at (17, 0.5) {};
		\node [style=none] (27) at (19, -0.5) {};
		\node [style=none] (28) at (19, 0.5) {};
		\node [style=none] (29) at (17, -0.5) {};
		\node [style=none] (30) at (10, 0.5) {};
		\node [style=none] (31) at (11, 0.5) {};
		\node [style=none] (32) at (11, -0.5) {};
		\node [style=none] (33) at (10, -0.5) {};
		\node [style=none] (34) at (10.5, 0) {$\bvdots$};
		\node [style=none] (35) at (25, 0.5) {};
		\node [style=none] (36) at (26, 0.5) {};
		\node [style=none] (37) at (26, -0.5) {};
		\node [style=none] (38) at (25, -0.5) {};
		\node [style=none] (39) at (2, 2) {$f: \sum\limits_{i} c_i  \to \sum\limits_{j} c_j$};
		\node [style=none] (40) at (10.5, 2) {$id_{\sum\limits_{i} c_i}: \sum\limits_{i} c_i  \to \sum\limits_{i} c_i$};
		\node [style=none] (41) at (18, 2) {$\sigma_{a,b}: a+b \to b+a$};
		\node [style=none] (42) at (25.5, 2) {$id_0 : 0\to 0$};
	\end{pgfonlayer}
	\begin{pgfonlayer}{edgelayer}
		\draw (0.center) to (1.center);
		\draw (1.center) to (3.center);
		\draw (3.center) to (2.center);
		\draw (2.center) to (0.center);
		\draw (4.center) to (9.center);
		\draw (5.center) to (8.center);
		\draw (6.center) to (7.center);
		\draw (17.center) to (18.center);
		\draw [in=-180, out=0, looseness=1.25] (26.center) to (27.center);
		\draw [in=0, out=-180, looseness=1.25] (28.center) to (29.center);
		\draw (30.center) to (31.center);
		\draw (32.center) to (33.center);
		\draw [style=pointille] (38.center)
			 to (37.center)
			 to (36.center)
			 to (35.center)
			 to cycle;
	\end{pgfonlayer}
\end{tikzpicture}

%% file: figures/stringdiagcomp.tikz
\begin{tikzpicture}
	\begin{pgfonlayer}{nodelayer}
		\node [style=none] (0) at (1, -0.5) {};
		\node [style=none] (1) at (1, -2.5) {};
		\node [style=none] (2) at (3, -0.5) {};
		\node [style=none] (3) at (3, -2.5) {};
		\node [style=none] (4) at (0, -1) {};
		\node [style=none] (5) at (0, -2) {};
		\node [style=none] (6) at (3, -1) {};
		\node [style=none] (7) at (4, -1) {};
		\node [style=none] (8) at (1, -2) {};
		\node [style=none] (9) at (1, -1) {};
		\node [style=none] (10) at (2, -1.5) {$f$};
		\node [style=none] (17) at (4, -2) {};
		\node [style=none] (18) at (3, -2) {};
		\node [style=none] (24) at (0.5, -1.5) {$\bvdots$};
		\node [style=none] (25) at (3.5, -1.5) {$\bvdots$};
		\node [style=none] (39) at (3.5, 2) {$g \circ f$};
		\node [style=none] (43) at (4, -0.5) {};
		\node [style=none] (44) at (4, -2.5) {};
		\node [style=none] (45) at (6, -0.5) {};
		\node [style=none] (46) at (6, -2.5) {};
		\node [style=none] (47) at (6, -1) {};
		\node [style=none] (48) at (7, -1) {};
		\node [style=none] (51) at (5, -1.5) {$g$};
		\node [style=none] (52) at (7, -2) {};
		\node [style=none] (53) at (6, -2) {};
		\node [style=none] (54) at (6.5, -1.5) {$\bvdots$};
		\node [style=none] (55) at (14, 1) {};
		\node [style=none] (56) at (14, -1) {};
		\node [style=none] (57) at (16, 1) {};
		\node [style=none] (58) at (16, -1) {};
		\node [style=none] (59) at (13, 0.5) {};
		\node [style=none] (60) at (13, -0.5) {};
		\node [style=none] (61) at (16, 0.5) {};
		\node [style=none] (62) at (17, 0.5) {};
		\node [style=none] (63) at (14, -0.5) {};
		\node [style=none] (64) at (14, 0.5) {};
		\node [style=none] (65) at (15, 0) {$f$};
		\node [style=none] (66) at (17, -0.5) {};
		\node [style=none] (67) at (16, -0.5) {};
		\node [style=none] (68) at (13.5, 0) {$\bvdots$};
		\node [style=none] (69) at (16.5, 0) {$\bvdots$};
		\node [style=none] (82) at (13, -2.5) {};
		\node [style=none] (83) at (14, -2.5) {};
		\node [style=none] (84) at (14, -3.5) {};
		\node [style=none] (85) at (13, -3.5) {};
		\node [style=none] (86) at (13.5, -3) {$\bvdots$};
		\node [style=none] (87) at (14, -2) {};
		\node [style=none] (88) at (14, -4) {};
		\node [style=none] (89) at (16, -2) {};
		\node [style=none] (90) at (16, -4) {};
		\node [style=none] (91) at (16, -2.5) {};
		\node [style=none] (92) at (17, -2.5) {};
		\node [style=none] (93) at (15, -3) {$g$};
		\node [style=none] (94) at (17, -3.5) {};
		\node [style=none] (95) at (16, -3.5) {};
		\node [style=none] (96) at (16.5, -3) {$\bvdots$};
		\node [style=none] (97) at (15, 2) {$f \otimes g $};
	\end{pgfonlayer}
	\begin{pgfonlayer}{edgelayer}
		\draw (0.center) to (1.center);
		\draw (1.center) to (3.center);
		\draw (3.center) to (2.center);
		\draw (2.center) to (0.center);
		\draw (4.center) to (9.center);
		\draw (5.center) to (8.center);
		\draw (6.center) to (7.center);
		\draw (17.center) to (18.center);
		\draw (43.center) to (44.center);
		\draw (44.center) to (46.center);
		\draw (46.center) to (45.center);
		\draw (45.center) to (43.center);
		\draw (47.center) to (48.center);
		\draw (52.center) to (53.center);
		\draw (55.center) to (56.center);
		\draw (56.center) to (58.center);
		\draw (58.center) to (57.center);
		\draw (57.center) to (55.center);
		\draw (59.center) to (64.center);
		\draw (60.center) to (63.center);
		\draw (61.center) to (62.center);
		\draw (66.center) to (67.center);
		\draw (82.center) to (83.center);
		\draw (84.center) to (85.center);
		\draw (87.center) to (88.center);
		\draw (88.center) to (90.center);
		\draw (90.center) to (89.center);
		\draw (89.center) to (87.center);
		\draw (91.center) to (92.center);
		\draw (94.center) to (95.center);
	\end{pgfonlayer}
\end{tikzpicture}

%% file: figures/stringdiageq.tikz
\begin{tikzpicture}
	\begin{pgfonlayer}{nodelayer}
		\node [style=none] (0) at (1, 4.5) {};
		\node [style=none] (1) at (1, 2.5) {};
		\node [style=none] (2) at (3, 4.5) {};
		\node [style=none] (3) at (3, 2.5) {};
		\node [style=none] (4) at (0, 4) {};
		\node [style=none] (5) at (0, 3) {};
		\node [style=none] (6) at (3, 4) {};
		\node [style=none] (7) at (4, 4) {};
		\node [style=none] (8) at (1, 3) {};
		\node [style=none] (9) at (1, 4) {};
		\node [style=none] (10) at (2, 3.5) {$f$};
		\node [style=none] (17) at (4, 3) {};
		\node [style=none] (18) at (3, 3) {};
		\node [style=none] (24) at (0.5, 3.5) {$\bvdots$};
		\node [style=none] (25) at (3.5, 3.5) {$\bvdots$};
		\node [style=none] (26) at (4, 2) {};
		\node [style=none] (27) at (6, 4) {};
		\node [style=none] (28) at (6, 3) {};
		\node [style=none] (29) at (4, 4) {};
		\node [style=none] (30) at (0, 2) {};
		\node [style=none] (31) at (4, 2) {};
		\node [style=none] (32) at (7, 4) {};
		\node [style=none] (33) at (6, 4) {};
		\node [style=none] (41) at (8, 3) {$=$};
		\node [style=none] (43) at (6, 3) {};
		\node [style=none] (44) at (7, 3) {};
		\node [style=none] (45) at (7, 2) {};
		\node [style=none] (46) at (6, 2) {};
		\node [style=none] (47) at (6.5, 2.5) {$\bvdots$};
		\node [style=none] (51) at (6, 2) {};
		\node [style=none] (52) at (4, 3) {};
		\node [style=none] (53) at (15, 1.5) {};
		\node [style=none] (54) at (15, 3.5) {};
		\node [style=none] (55) at (13, 1.5) {};
		\node [style=none] (56) at (13, 3.5) {};
		\node [style=none] (57) at (16, 2) {};
		\node [style=none] (58) at (16, 3) {};
		\node [style=none] (59) at (13, 2) {};
		\node [style=none] (60) at (12, 2) {};
		\node [style=none] (61) at (15, 3) {};
		\node [style=none] (62) at (15, 2) {};
		\node [style=none] (63) at (14, 2.5) {$f$};
		\node [style=none] (64) at (12, 3) {};
		\node [style=none] (65) at (13, 3) {};
		\node [style=none] (66) at (15.5, 2.5) {$\bvdots$};
		\node [style=none] (67) at (12.5, 2.5) {$\bvdots$};
		\node [style=none] (68) at (12, 4) {};
		\node [style=none] (69) at (10, 2) {};
		\node [style=none] (70) at (10, 3) {};
		\node [style=none] (71) at (12, 2) {};
		\node [style=none] (72) at (16, 4) {};
		\node [style=none] (73) at (12, 4) {};
		\node [style=none] (74) at (9, 2) {};
		\node [style=none] (75) at (10, 2) {};
		\node [style=none] (76) at (10, 3) {};
		\node [style=none] (77) at (9, 3) {};
		\node [style=none] (78) at (9, 4) {};
		\node [style=none] (79) at (10, 4) {};
		\node [style=none] (80) at (9.5, 3.5) {$\bvdots$};
		\node [style=none] (81) at (10, 4) {};
		\node [style=none] (82) at (12, 3) {};
		\node [style=none] (83) at (21, 4) {};
		\node [style=none] (84) at (21, 2) {};
		\node [style=none] (85) at (23, 4) {};
		\node [style=none] (86) at (23, 2) {};
		\node [style=none] (87) at (20, 3.5) {};
		\node [style=none] (88) at (20, 2.5) {};
		\node [style=none] (89) at (23, 3.5) {};
		\node [style=none] (90) at (24, 3.5) {};
		\node [style=none] (91) at (21, 2.5) {};
		\node [style=none] (92) at (21, 3.5) {};
		\node [style=none] (93) at (22, 3) {$f$};
		\node [style=none] (94) at (24, 2.5) {};
		\node [style=none] (95) at (23, 2.5) {};
		\node [style=none] (96) at (20.5, 3) {$\bvdots$};
		\node [style=none] (97) at (23.5, 3) {$\bvdots$};
		\node [style=none] (100) at (26, 3.5) {};
		\node [style=none] (101) at (24, 3.5) {};
		\node [style=none] (106) at (28, 3) {$=$};
		\node [style=none] (107) at (26, 3.5) {};
		\node [style=none] (108) at (27, 3.5) {};
		\node [style=none] (109) at (27, 2.5) {};
		\node [style=none] (110) at (26, 2.5) {};
		\node [style=none] (111) at (26.5, 3) {$\bvdots$};
		\node [style=none] (112) at (26, 2.5) {};
		\node [style=none] (113) at (24, 2.5) {};
		\node [style=none] (114) at (35, 2) {};
		\node [style=none] (115) at (35, 4) {};
		\node [style=none] (116) at (33, 2) {};
		\node [style=none] (117) at (33, 4) {};
		\node [style=none] (118) at (36, 2.5) {};
		\node [style=none] (119) at (36, 3.5) {};
		\node [style=none] (120) at (33, 2.5) {};
		\node [style=none] (121) at (32, 2.5) {};
		\node [style=none] (122) at (35, 3.5) {};
		\node [style=none] (123) at (35, 2.5) {};
		\node [style=none] (124) at (34, 3) {$f$};
		\node [style=none] (125) at (32, 3.5) {};
		\node [style=none] (126) at (33, 3.5) {};
		\node [style=none] (127) at (35.5, 3) {$\bvdots$};
		\node [style=none] (128) at (32.5, 3) {$\bvdots$};
		\node [style=none] (131) at (30, 2.5) {};
		\node [style=none] (132) at (32, 2.5) {};
		\node [style=none] (137) at (30, 2.5) {};
		\node [style=none] (138) at (29, 2.5) {};
		\node [style=none] (139) at (29, 3.5) {};
		\node [style=none] (140) at (30, 3.5) {};
		\node [style=none] (141) at (29.5, 3) {$\bvdots$};
		\node [style=none] (142) at (30, 3.5) {};
		\node [style=none] (143) at (32, 3.5) {};
		\node [style=none] (161) at (5, -3) {$=$};
		\node [style=none] (169) at (9, -4) {};
		\node [style=none] (170) at (9, -2) {};
		\node [style=none] (171) at (7, -4) {};
		\node [style=none] (172) at (7, -2) {};
		\node [style=none] (173) at (10, -3.5) {};
		\node [style=none] (174) at (10, -2.5) {};
		\node [style=none] (175) at (7, -3.5) {};
		\node [style=none] (176) at (6, -3.5) {};
		\node [style=none] (177) at (9, -2.5) {};
		\node [style=none] (178) at (9, -3.5) {};
		\node [style=none] (179) at (8, -3) {$f$};
		\node [style=none] (180) at (6, -2.5) {};
		\node [style=none] (181) at (7, -2.5) {};
		\node [style=none] (182) at (9.5, -3) {$\bvdots$};
		\node [style=none] (183) at (6.5, -3) {$\bvdots$};
		\node [style=none] (185) at (6, -3.5) {};
		\node [style=none] (192) at (6, -2.5) {};
		\node [style=none] (193) at (13, -3.5) {};
		\node [style=none] (194) at (13, -1.5) {};
		\node [style=none] (195) at (15, -3.5) {};
		\node [style=none] (196) at (15, -1.5) {};
		\node [style=none] (197) at (12, -3) {};
		\node [style=none] (198) at (12, -2) {};
		\node [style=none] (199) at (15, -3) {};
		\node [style=none] (200) at (16, -3) {};
		\node [style=none] (201) at (13, -2) {};
		\node [style=none] (202) at (13, -3) {};
		\node [style=none] (203) at (14, -2.5) {$f$};
		\node [style=none] (204) at (16, -2) {};
		\node [style=none] (205) at (15, -2) {};
		\node [style=none] (206) at (12.5, -2.5) {$\bvdots$};
		\node [style=none] (207) at (15.5, -2.5) {$\bvdots$};
		\node [style=none] (208) at (16, -3) {};
		\node [style=none] (209) at (16, -2) {};
		\node [style=none] (210) at (13.5, -5) {};
		\node [style=none] (211) at (14.5, -5) {};
		\node [style=none] (212) at (14.5, -4) {};
		\node [style=none] (213) at (13.5, -4) {};
		\node [style=none] (214) at (1, -3) {};
		\node [style=none] (215) at (1, -5) {};
		\node [style=none] (216) at (3, -3) {};
		\node [style=none] (217) at (3, -5) {};
		\node [style=none] (218) at (0, -3.5) {};
		\node [style=none] (219) at (0, -4.5) {};
		\node [style=none] (220) at (3, -3.5) {};
		\node [style=none] (221) at (4, -3.5) {};
		\node [style=none] (222) at (1, -4.5) {};
		\node [style=none] (223) at (1, -3.5) {};
		\node [style=none] (224) at (2, -4) {$f$};
		\node [style=none] (225) at (4, -4.5) {};
		\node [style=none] (226) at (3, -4.5) {};
		\node [style=none] (227) at (0.5, -4) {$\bvdots$};
		\node [style=none] (228) at (3.5, -4) {$\bvdots$};
		\node [style=none] (229) at (4, -3.5) {};
		\node [style=none] (230) at (4, -4.5) {};
		\node [style=none] (231) at (1.5, -1.5) {};
		\node [style=none] (232) at (2.5, -1.5) {};
		\node [style=none] (233) at (2.5, -2.5) {};
		\node [style=none] (234) at (1.5, -2.5) {};
		\node [style=none] (235) at (11, -3) {$=$};
		\node [style=none] (236) at (25, 4.5) {};
		\node [style=none] (237) at (25, 1.5) {};
		\node [style=none] (238) at (31, 4.5) {};
		\node [style=none] (239) at (31, 1.5) {};
		\node [style=none] (240) at (22, -2) {};
		\node [style=none] (241) at (22, 0) {};
		\node [style=none] (242) at (20, -2) {};
		\node [style=none] (243) at (20, 0) {};
		\node [style=none] (244) at (23, -1.5) {};
		\node [style=none] (245) at (23, -0.5) {};
		\node [style=none] (246) at (20, -1.5) {};
		\node [style=none] (247) at (19, -1.5) {};
		\node [style=none] (248) at (22, -0.5) {};
		\node [style=none] (249) at (22, -1.5) {};
		\node [style=none] (250) at (21, -1) {$f$};
		\node [style=none] (251) at (19, -0.5) {};
		\node [style=none] (252) at (20, -0.5) {};
		\node [style=none] (253) at (22.5, -1) {$\bvdots$};
		\node [style=none] (254) at (19.5, -1) {$\bvdots$};
		\node [style=none] (255) at (22, -6) {};
		\node [style=none] (256) at (22, -4) {};
		\node [style=none] (257) at (20, -6) {};
		\node [style=none] (258) at (20, -4) {};
		\node [style=none] (259) at (23, -5.5) {};
		\node [style=none] (260) at (23, -4.5) {};
		\node [style=none] (261) at (20, -5.5) {};
		\node [style=none] (262) at (19, -5.5) {};
		\node [style=none] (263) at (22, -4.5) {};
		\node [style=none] (264) at (22, -5.5) {};
		\node [style=none] (265) at (21, -5) {$g$};
		\node [style=none] (266) at (19, -4.5) {};
		\node [style=none] (267) at (20, -4.5) {};
		\node [style=none] (268) at (22.5, -5) {$\bvdots$};
		\node [style=none] (269) at (19.5, -5) {$\bvdots$};
		\node [style=none] (270) at (26, -2) {};
		\node [style=none] (271) at (26, 0) {};
		\node [style=none] (272) at (24, -2) {};
		\node [style=none] (273) at (24, 0) {};
		\node [style=none] (274) at (27, -1.5) {};
		\node [style=none] (275) at (27, -0.5) {};
		\node [style=none] (276) at (24, -1.5) {};
		\node [style=none] (277) at (23, -1.5) {};
		\node [style=none] (278) at (26, -0.5) {};
		\node [style=none] (279) at (26, -1.5) {};
		\node [style=none] (280) at (25, -1) {$k$};
		\node [style=none] (281) at (23, -0.5) {};
		\node [style=none] (282) at (24, -0.5) {};
		\node [style=none] (283) at (26.5, -1) {$\bvdots$};
		\node [style=none] (284) at (23.5, -1) {$\bvdots$};
		\node [style=none] (285) at (26, -6) {};
		\node [style=none] (286) at (26, -4) {};
		\node [style=none] (287) at (24, -6) {};
		\node [style=none] (288) at (24, -4) {};
		\node [style=none] (289) at (27, -5.5) {};
		\node [style=none] (290) at (27, -4.5) {};
		\node [style=none] (291) at (24, -5.5) {};
		\node [style=none] (292) at (23, -5.5) {};
		\node [style=none] (293) at (26, -4.5) {};
		\node [style=none] (294) at (26, -5.5) {};
		\node [style=none] (295) at (25, -5) {$h$};
		\node [style=none] (296) at (23, -4.5) {};
		\node [style=none] (297) at (24, -4.5) {};
		\node [style=none] (298) at (26.5, -5) {$\bvdots$};
		\node [style=none] (299) at (23.5, -5) {$\bvdots$};
		\node [style=none] (300) at (32, -2) {};
		\node [style=none] (301) at (32, 0) {};
		\node [style=none] (302) at (30, -2) {};
		\node [style=none] (303) at (30, 0) {};
		\node [style=none] (304) at (33, -1.5) {};
		\node [style=none] (305) at (33, -0.5) {};
		\node [style=none] (306) at (30, -1.5) {};
		\node [style=none] (307) at (29, -1.5) {};
		\node [style=none] (308) at (32, -0.5) {};
		\node [style=none] (309) at (32, -1.5) {};
		\node [style=none] (310) at (31, -1) {$f$};
		\node [style=none] (311) at (29, -0.5) {};
		\node [style=none] (312) at (30, -0.5) {};
		\node [style=none] (313) at (32.5, -1) {$\bvdots$};
		\node [style=none] (314) at (29.5, -1) {$\bvdots$};
		\node [style=none] (315) at (32, -6) {};
		\node [style=none] (316) at (32, -4) {};
		\node [style=none] (317) at (30, -6) {};
		\node [style=none] (318) at (30, -4) {};
		\node [style=none] (319) at (33, -5.5) {};
		\node [style=none] (320) at (33, -4.5) {};
		\node [style=none] (321) at (30, -5.5) {};
		\node [style=none] (322) at (29, -5.5) {};
		\node [style=none] (323) at (32, -4.5) {};
		\node [style=none] (324) at (32, -5.5) {};
		\node [style=none] (325) at (31, -5) {$g$};
		\node [style=none] (326) at (29, -4.5) {};
		\node [style=none] (327) at (30, -4.5) {};
		\node [style=none] (328) at (32.5, -5) {$\bvdots$};
		\node [style=none] (329) at (29.5, -5) {$\bvdots$};
		\node [style=none] (330) at (36, -2) {};
		\node [style=none] (331) at (36, 0) {};
		\node [style=none] (332) at (34, -2) {};
		\node [style=none] (333) at (34, 0) {};
		\node [style=none] (334) at (37, -1.5) {};
		\node [style=none] (335) at (37, -0.5) {};
		\node [style=none] (336) at (34, -1.5) {};
		\node [style=none] (337) at (33, -1.5) {};
		\node [style=none] (338) at (36, -0.5) {};
		\node [style=none] (339) at (36, -1.5) {};
		\node [style=none] (340) at (35, -1) {$k$};
		\node [style=none] (341) at (33, -0.5) {};
		\node [style=none] (342) at (34, -0.5) {};
		\node [style=none] (343) at (36.5, -1) {$\bvdots$};
		\node [style=none] (344) at (33.5, -1) {$\bvdots$};
		\node [style=none] (345) at (36, -6) {};
		\node [style=none] (346) at (36, -4) {};
		\node [style=none] (347) at (34, -6) {};
		\node [style=none] (348) at (34, -4) {};
		\node [style=none] (349) at (37, -5.5) {};
		\node [style=none] (350) at (37, -4.5) {};
		\node [style=none] (351) at (34, -5.5) {};
		\node [style=none] (352) at (33, -5.5) {};
		\node [style=none] (353) at (36, -4.5) {};
		\node [style=none] (354) at (36, -5.5) {};
		\node [style=none] (355) at (35, -5) {$h$};
		\node [style=none] (356) at (33, -4.5) {};
		\node [style=none] (357) at (34, -4.5) {};
		\node [style=none] (358) at (36.5, -5) {$\bvdots$};
		\node [style=none] (359) at (33.5, -5) {$\bvdots$};
		\node [style=none] (360) at (33, 0.5) {};
		\node [style=none] (361) at (33, -6.5) {};
		\node [style=none] (362) at (27, -3) {};
		\node [style=none] (363) at (19, -3) {};
		\node [style=none] (364) at (28, -3) {$=$};
	\end{pgfonlayer}
	\begin{pgfonlayer}{edgelayer}
		\draw (0.center) to (1.center);
		\draw (1.center) to (3.center);
		\draw (3.center) to (2.center);
		\draw (2.center) to (0.center);
		\draw (4.center) to (9.center);
		\draw (5.center) to (8.center);
		\draw (6.center) to (7.center);
		\draw (17.center) to (18.center);
		\draw [in=-180, out=0, looseness=1.25] (26.center) to (27.center);
		\draw [in=0, out=-180, looseness=1.25] (28.center) to (29.center);
		\draw (30.center) to (31.center);
		\draw (32.center) to (33.center);
		\draw (43.center) to (44.center);
		\draw (45.center) to (46.center);
		\draw [in=0, out=-180, looseness=1.25] (51.center) to (52.center);
		\draw (53.center) to (54.center);
		\draw (54.center) to (56.center);
		\draw (56.center) to (55.center);
		\draw (55.center) to (53.center);
		\draw (57.center) to (62.center);
		\draw (58.center) to (61.center);
		\draw (59.center) to (60.center);
		\draw (64.center) to (65.center);
		\draw [in=0, out=-180, looseness=1.25] (68.center) to (69.center);
		\draw [in=-180, out=0, looseness=1.25] (70.center) to (71.center);
		\draw (72.center) to (73.center);
		\draw (74.center) to (75.center);
		\draw (76.center) to (77.center);
		\draw (78.center) to (79.center);
		\draw [in=-180, out=0, looseness=1.25] (81.center) to (82.center);
		\draw (83.center) to (84.center);
		\draw (84.center) to (86.center);
		\draw (86.center) to (85.center);
		\draw (85.center) to (83.center);
		\draw (87.center) to (92.center);
		\draw (88.center) to (91.center);
		\draw (89.center) to (90.center);
		\draw (94.center) to (95.center);
		\draw [in=0, out=-180, looseness=1.25] (100.center) to (101.center);
		\draw (107.center) to (108.center);
		\draw (109.center) to (110.center);
		\draw [in=0, out=-180, looseness=1.25] (112.center) to (113.center);
		\draw (114.center) to (115.center);
		\draw (115.center) to (117.center);
		\draw (117.center) to (116.center);
		\draw (116.center) to (114.center);
		\draw (118.center) to (123.center);
		\draw (119.center) to (122.center);
		\draw (120.center) to (121.center);
		\draw (125.center) to (126.center);
		\draw [in=-180, out=0, looseness=1.25] (131.center) to (132.center);
		\draw (137.center) to (138.center);
		\draw (139.center) to (140.center);
		\draw [in=-180, out=0, looseness=1.25] (142.center) to (143.center);
		\draw (169.center) to (170.center);
		\draw (170.center) to (172.center);
		\draw (172.center) to (171.center);
		\draw (171.center) to (169.center);
		\draw (173.center) to (178.center);
		\draw (174.center) to (177.center);
		\draw (175.center) to (176.center);
		\draw (180.center) to (181.center);
		\draw (193.center) to (194.center);
		\draw (194.center) to (196.center);
		\draw (196.center) to (195.center);
		\draw (195.center) to (193.center);
		\draw (197.center) to (202.center);
		\draw (198.center) to (201.center);
		\draw (199.center) to (200.center);
		\draw (204.center) to (205.center);
		\draw [style=pointille] (213.center) to (212.center);
		\draw [style=pointille] (212.center) to (211.center);
		\draw [style=pointille] (211.center) to (210.center);
		\draw [style=pointille] (210.center) to (213.center);
		\draw (214.center) to (215.center);
		\draw (215.center) to (217.center);
		\draw (217.center) to (216.center);
		\draw (216.center) to (214.center);
		\draw (218.center) to (223.center);
		\draw (219.center) to (222.center);
		\draw (220.center) to (221.center);
		\draw (225.center) to (226.center);
		\draw [style=pointille] (234.center) to (233.center);
		\draw [style=pointille] (233.center) to (232.center);
		\draw [style=pointille] (232.center) to (231.center);
		\draw [style=pointille] (231.center) to (234.center);
		\draw [style=pointille] (236.center) to (237.center);
		\draw [style=pointille] (239.center) to (238.center);
		\draw (240.center) to (241.center);
		\draw (241.center) to (243.center);
		\draw (243.center) to (242.center);
		\draw (242.center) to (240.center);
		\draw (244.center) to (249.center);
		\draw (245.center) to (248.center);
		\draw (246.center) to (247.center);
		\draw (251.center) to (252.center);
		\draw (255.center) to (256.center);
		\draw (256.center) to (258.center);
		\draw (258.center) to (257.center);
		\draw (257.center) to (255.center);
		\draw (259.center) to (264.center);
		\draw (260.center) to (263.center);
		\draw (261.center) to (262.center);
		\draw (266.center) to (267.center);
		\draw (270.center) to (271.center);
		\draw (271.center) to (273.center);
		\draw (273.center) to (272.center);
		\draw (272.center) to (270.center);
		\draw (274.center) to (279.center);
		\draw (275.center) to (278.center);
		\draw (276.center) to (277.center);
		\draw (281.center) to (282.center);
		\draw (285.center) to (286.center);
		\draw (286.center) to (288.center);
		\draw (288.center) to (287.center);
		\draw (287.center) to (285.center);
		\draw (289.center) to (294.center);
		\draw (290.center) to (293.center);
		\draw (291.center) to (292.center);
		\draw (296.center) to (297.center);
		\draw (300.center) to (301.center);
		\draw (301.center) to (303.center);
		\draw (303.center) to (302.center);
		\draw (302.center) to (300.center);
		\draw (304.center) to (309.center);
		\draw (305.center) to (308.center);
		\draw (306.center) to (307.center);
		\draw (311.center) to (312.center);
		\draw (315.center) to (316.center);
		\draw (316.center) to (318.center);
		\draw (318.center) to (317.center);
		\draw (317.center) to (315.center);
		\draw (319.center) to (324.center);
		\draw (320.center) to (323.center);
		\draw (321.center) to (322.center);
		\draw (326.center) to (327.center);
		\draw (330.center) to (331.center);
		\draw (331.center) to (333.center);
		\draw (333.center) to (332.center);
		\draw (332.center) to (330.center);
		\draw (334.center) to (339.center);
		\draw (335.center) to (338.center);
		\draw (336.center) to (337.center);
		\draw (341.center) to (342.center);
		\draw (345.center) to (346.center);
		\draw (346.center) to (348.center);
		\draw (348.center) to (347.center);
		\draw (347.center) to (345.center);
		\draw (349.center) to (354.center);
		\draw (350.center) to (353.center);
		\draw (351.center) to (352.center);
		\draw (356.center) to (357.center);
		\draw [style=pointille] (361.center) to (360.center);
		\draw [style=pointille] (363.center) to (362.center);
	\end{pgfonlayer}
\end{tikzpicture}

%% file: figures/bureaucracydef.tikz
\begin{tikzpicture}
	\begin{pgfonlayer}{nodelayer}
		\node [style=none] (30) at (0, 0.5) {};
		\node [style=none] (31) at (1, 0.5) {};
		\node [style=none] (32) at (1, -0.5) {};
		\node [style=none] (33) at (0, -0.5) {};
		\node [style=none] (34) at (0.5, 0) {$\bvdots$};
		\node [style=none] (44) at (2, 0) {};
		\node [style=none] (46) at (3, 0.5) {};
		\node [style=none] (47) at (4, 0.5) {};
		\node [style=none] (48) at (4, -0.5) {};
		\node [style=none] (49) at (3, -0.5) {};
		\node [style=none] (50) at (3.5, 0) {$\bvdots$};
		\node [style=new style 0] (51) at (2, 0) {};
	\end{pgfonlayer}
	\begin{pgfonlayer}{edgelayer}
		\draw (30.center) to (31.center);
		\draw (32.center) to (33.center);
		\draw (46.center) to (47.center);
		\draw (48.center) to (49.center);
		\draw [in=0, out=135, looseness=1.25] (44.center) to (31.center);
		\draw [in=-135, out=0, looseness=1.25] (32.center) to (44.center);
		\draw [in=180, out=-45, looseness=1.25] (44.center) to (49.center);
		\draw [in=45, out=-180, looseness=1.25] (46.center) to (44.center);
	\end{pgfonlayer}
\end{tikzpicture}

%% file: figures/bureaucracycalc.tikz
\begin{tikzpicture}
	\begin{pgfonlayer}{nodelayer}
		\node [style=none] (30) at (0, 0.5) {};
		\node [style=none] (31) at (1, 0.5) {};
		\node [style=none] (32) at (1, -0.5) {};
		\node [style=none] (33) at (0, -0.5) {};
		\node [style=none] (34) at (0.5, 0) {$\bvdots$};
		\node [style=none] (39) at (11, 0) {$=$};
		\node [style=none] (44) at (2, 0) {};
		\node [style=none] (46) at (3, 0.5) {};
		\node [style=none] (47) at (4, 0.5) {};
		\node [style=none] (48) at (4, -0.5) {};
		\node [style=none] (49) at (3, -0.5) {};
		\node [style=none] (50) at (3.5, 0) {$\bvdots$};
		\node [style=new style 0] (51) at (2, 0) {};
		\node [style=none] (53) at (4, 0.5) {};
		\node [style=none] (54) at (4, -0.5) {};
		\node [style=none] (57) at (5, 0) {};
		\node [style=none] (58) at (6, 0.5) {};
		\node [style=none] (61) at (6, -0.5) {};
		\node [style=new style 0] (63) at (5, 0) {};
		\node [style=none] (64) at (6, 0.5) {};
		\node [style=none] (65) at (7, 0.5) {};
		\node [style=none] (66) at (7, -0.5) {};
		\node [style=none] (67) at (6, -0.5) {};
		\node [style=none] (68) at (6.5, 0) {$\bvdots$};
		\node [style=none] (69) at (8, 0) {};
		\node [style=none] (70) at (9, 0.5) {};
		\node [style=none] (71) at (10, 0.5) {};
		\node [style=none] (72) at (10, -0.5) {};
		\node [style=none] (73) at (9, -0.5) {};
		\node [style=none] (74) at (9.5, 0) {$\bvdots$};
		\node [style=new style 0] (75) at (8, 0) {};
		\node [style=none] (76) at (12, 0.5) {};
		\node [style=none] (77) at (13, 0.5) {};
		\node [style=none] (78) at (13, -0.5) {};
		\node [style=none] (79) at (12, -0.5) {};
		\node [style=none] (80) at (12.5, 0) {$\bvdots$};
		\node [style=none] (81) at (14, 0) {};
		\node [style=none] (82) at (15, 0.5) {};
		\node [style=none] (83) at (16, 0.5) {};
		\node [style=none] (84) at (16, -0.5) {};
		\node [style=none] (85) at (15, -0.5) {};
		\node [style=none] (86) at (15.5, 0) {$\bvdots$};
		\node [style=new style 0] (87) at (14, 0) {};
		\node [style=none] (88) at (6, -1.5) {};
		\node [style=none] (89) at (7, -1.5) {};
		\node [style=none] (90) at (7, -2.5) {};
		\node [style=none] (91) at (6, -2.5) {};
		\node [style=none] (92) at (6.5, -2) {$\bvdots$};
		\node [style=none] (93) at (7, -1.5) {};
		\node [style=none] (94) at (7, -2.5) {};
		\node [style=none] (95) at (3, 2.5) {};
		\node [style=none] (96) at (4, 2.5) {};
		\node [style=none] (97) at (4, 1.5) {};
		\node [style=none] (98) at (3, 1.5) {};
		\node [style=none] (99) at (3.5, 2) {$\bvdots$};
		\node [style=none] (100) at (4, 2.5) {};
		\node [style=none] (101) at (4, 1.5) {};
		\node [style=none] (102) at (9, 2.5) {};
		\node [style=none] (103) at (10, 2.5) {};
		\node [style=none] (104) at (10, 1.5) {};
		\node [style=none] (105) at (9, 1.5) {};
		\node [style=none] (106) at (9.5, 2) {$\bvdots$};
		\node [style=none] (107) at (10, 2.5) {};
		\node [style=none] (108) at (10, 1.5) {};
		\node [style=none] (109) at (0, -1.5) {};
		\node [style=none] (110) at (1, -1.5) {};
		\node [style=none] (111) at (1, -2.5) {};
		\node [style=none] (112) at (0, -2.5) {};
		\node [style=none] (113) at (0.5, -2) {$\bvdots$};
		\node [style=none] (114) at (1, -1.5) {};
		\node [style=none] (115) at (1, -2.5) {};
	\end{pgfonlayer}
	\begin{pgfonlayer}{edgelayer}
		\draw (30.center) to (31.center);
		\draw (32.center) to (33.center);
		\draw (46.center) to (47.center);
		\draw (48.center) to (49.center);
		\draw [in=0, out=135, looseness=1.25] (44.center) to (31.center);
		\draw [in=-135, out=0, looseness=1.25] (32.center) to (44.center);
		\draw [in=180, out=-45, looseness=1.25] (44.center) to (49.center);
		\draw [in=45, out=-180, looseness=1.25] (46.center) to (44.center);
		\draw [in=0, out=135, looseness=1.25] (57.center) to (53.center);
		\draw [in=-135, out=0, looseness=1.25] (54.center) to (57.center);
		\draw [in=180, out=-45, looseness=1.25] (57.center) to (61.center);
		\draw [in=45, out=-180, looseness=1.25] (58.center) to (57.center);
		\draw (64.center) to (65.center);
		\draw (66.center) to (67.center);
		\draw (70.center) to (71.center);
		\draw (72.center) to (73.center);
		\draw [in=0, out=135, looseness=1.25] (69.center) to (65.center);
		\draw [in=-135, out=0, looseness=1.25] (66.center) to (69.center);
		\draw [in=180, out=-45, looseness=1.25] (69.center) to (73.center);
		\draw [in=45, out=-180, looseness=1.25] (70.center) to (69.center);
		\draw (76.center) to (77.center);
		\draw (78.center) to (79.center);
		\draw (82.center) to (83.center);
		\draw (84.center) to (85.center);
		\draw [in=0, out=135, looseness=1.25] (81.center) to (77.center);
		\draw [in=-135, out=0, looseness=1.25] (78.center) to (81.center);
		\draw [in=180, out=-45, looseness=1.25] (81.center) to (85.center);
		\draw [in=45, out=-180, looseness=1.25] (82.center) to (81.center);
		\draw (88.center) to (89.center);
		\draw (90.center) to (91.center);
		\draw (95.center) to (96.center);
		\draw (97.center) to (98.center);
		\draw [in=180, out=75, looseness=0.75] (51) to (95.center);
		\draw [in=-180, out=60] (51) to (98.center);
		\draw [in=240, out=0] (93.center) to (75);
		\draw [in=255, out=0, looseness=0.75] (94.center) to (75);
		\draw (102.center) to (103.center);
		\draw (104.center) to (105.center);
		\draw (109.center) to (110.center);
		\draw (111.center) to (112.center);
		\draw (100.center) to (102.center);
		\draw (101.center) to (105.center);
		\draw (88.center) to (114.center);
		\draw (115.center) to (91.center);
	\end{pgfonlayer}
\end{tikzpicture}

%% file: figures/divandgathdef.tikz
\begin{tikzpicture}
	\begin{pgfonlayer}{nodelayer}
		\node [style=none] (90) at (21, 0) {};
		\node [style=none] (91) at (20, 0) {};
		\node [style=none] (93) at (1, 0) {};
		\node [style=none] (94) at (0, 0) {};
		\node [style=none] (106) at (8, 0.5) {};
		\node [style=none] (107) at (8, -0.5) {};
		\node [style=none] (108) at (7, 0) {};
		\node [style=none] (109) at (6, 0) {};
		\node [style=none] (130) at (13, 0.5) {};
		\node [style=none] (131) at (13, -0.5) {};
		\node [style=none] (132) at (14, 0) {};
		\node [style=none] (133) at (15, 0) {};
		\node [style=divide] (137) at (7, 0) {};
		\node [style=divide] (142) at (0, 0) {};
		\node [style=gather] (146) at (14, 0) {};
		\node [style=gather] (147) at (21, 0) {};
		\node [style=none] (148) at (20.5, 1) {$[I]\to 0$};
		\node [style=none] (149) at (0.5, 1) {$0 \to [I]$};
		\node [style=none] (150) at (14, 1.5) {$[a]+[b] \to  [a\otimes b]$};
		\node [style=none] (151) at (7, 1.5) {$[a\otimes b] \to   [a]+[b] $};
	\end{pgfonlayer}
	\begin{pgfonlayer}{edgelayer}
		\draw (90.center) to (91.center);
		\draw (93.center) to (94.center);
		\draw [in=-180, out=45, looseness=1.25] (108.center) to (106.center);
		\draw [in=-45, out=180, looseness=1.25] (107.center) to (108.center);
		\draw [in=-180, out=0, looseness=1.25] (109.center) to (108.center);
		\draw [in=0, out=135, looseness=1.25] (132.center) to (130.center);
		\draw [in=-135, out=0, looseness=1.25] (131.center) to (132.center);
		\draw [in=0, out=180, looseness=1.25] (133.center) to (132.center);
	\end{pgfonlayer}
\end{tikzpicture}

%% file: figures/divandgatheq.tikz
\begin{tikzpicture}
	\begin{pgfonlayer}{nodelayer}
		\node [style=none] (51) at (12.5, 2) {$=$};
		\node [style=none] (68) at (9.75, 2.5) {};
		\node [style=none] (69) at (9.75, 1.5) {};
		\node [style=none] (70) at (8.75, 2) {};
		\node [style=none] (71) at (7.75, 2) {};
		\node [style=none] (72) at (9.75, 2.5) {};
		\node [style=none] (73) at (9.75, 1.5) {};
		\node [style=none] (74) at (10.75, 2) {};
		\node [style=none] (75) at (11.75, 2) {};
		\node [style=none] (76) at (11.75, -1.5) {};
		\node [style=none] (77) at (11.75, -2.5) {};
		\node [style=none] (78) at (10.75, -2) {};
		\node [style=none] (79) at (9.75, -2) {};
		\node [style=none] (80) at (8.75, -1.5) {};
		\node [style=none] (81) at (8.75, -2.5) {};
		\node [style=none] (82) at (9.75, -2) {};
		\node [style=none] (83) at (12.5, -2) {$=$};
		\node [style=none] (84) at (13.25, 2) {};
		\node [style=none] (85) at (15.25, 2) {};
		\node [style=none] (86) at (13.25, -1.5) {};
		\node [style=none] (87) at (15.25, -1.5) {};
		\node [style=none] (88) at (15.25, -2.5) {};
		\node [style=none] (89) at (13.25, -2.5) {};
		\node [style=none] (90) at (18.75, 2) {};
		\node [style=none] (91) at (17.75, 2) {};
		\node [style=none] (92) at (21.5, 2) {$=$};
		\node [style=none] (93) at (20.75, 2) {};
		\node [style=none] (94) at (19.75, 2) {};
		\node [style=none] (95) at (24.25, 2) {};
		\node [style=none] (96) at (22.25, 2) {};
		\node [style=none] (99) at (21.5, -2) {$=$};
		\node [style=none] (100) at (20.75, -2) {};
		\node [style=none] (101) at (19.75, -2) {};
		\node [style=none] (102) at (23.25, -1.5) {};
		\node [style=none] (103) at (22.25, -1.5) {};
		\node [style=none] (104) at (23.25, -2.5) {};
		\node [style=none] (105) at (22.25, -2.5) {};
		\node [style=none] (106) at (2, 2.5) {};
		\node [style=none] (107) at (2, 1.5) {};
		\node [style=none] (108) at (1, 2) {};
		\node [style=none] (109) at (0, 2) {};
		\node [style=none] (111) at (2.75, 2) {$=$};
		\node [style=none] (114) at (5.5, 2) {};
		\node [style=none] (115) at (3.5, 2) {};
		\node [style=none] (116) at (2, -1.5) {};
		\node [style=none] (117) at (2, -2.5) {};
		\node [style=none] (118) at (1, -2) {};
		\node [style=none] (119) at (0, -2) {};
		\node [style=none] (120) at (2.75, -2) {$=$};
		\node [style=none] (121) at (5.5, -2) {};
		\node [style=none] (122) at (3.5, -2) {};
		\node [style=none] (123) at (29.75, 2.25) {};
		\node [style=none] (124) at (29.75, 1.25) {};
		\node [style=none] (125) at (30.75, 1.75) {};
		\node [style=none] (126) at (31.75, 1.75) {};
		\node [style=none] (127) at (29, 1.75) {$=$};
		\node [style=none] (128) at (26.25, 1.75) {};
		\node [style=none] (129) at (28.25, 1.75) {};
		\node [style=none] (130) at (29.75, -1.75) {};
		\node [style=none] (131) at (29.75, -2.75) {};
		\node [style=none] (132) at (30.75, -2.25) {};
		\node [style=none] (133) at (31.75, -2.25) {};
		\node [style=none] (134) at (29, -2.25) {$=$};
		\node [style=none] (135) at (26.25, -2.25) {};
		\node [style=none] (136) at (28.25, -2.25) {};
		\node [style=divide] (137) at (1, 2) {};
		\node [style=divide] (138) at (1, -2) {};
		\node [style=divide] (139) at (8.75, 2) {};
		\node [style=divide] (140) at (10.75, -2) {};
		\node [style=divide] (141) at (19.75, -2) {};
		\node [style=divide] (142) at (19.75, 2) {};
		\node [style=divide] (143) at (29.75, 2.25) {};
		\node [style=divide] (144) at (29.75, -2.75) {};
		\node [style=gather] (145) at (30.75, 1.75) {};
		\node [style=gather] (146) at (30.75, -2.25) {};
		\node [style=gather] (147) at (18.75, 2) {};
		\node [style=gather] (148) at (20.75, -2) {};
		\node [style=gather] (149) at (10.75, 2) {};
		\node [style=gather] (150) at (9.75, -2) {};
		\node [style=gather] (151) at (2, 2.5) {};
		\node [style=gather] (152) at (2, -2.5) {};
	\end{pgfonlayer}
	\begin{pgfonlayer}{edgelayer}
		\draw [in=-180, out=45, looseness=1.25] (70.center) to (68.center);
		\draw [in=-45, out=180, looseness=1.25] (69.center) to (70.center);
		\draw [in=-180, out=0, looseness=1.25] (71.center) to (70.center);
		\draw [in=0, out=135, looseness=1.25] (74.center) to (72.center);
		\draw [in=-135, out=0, looseness=1.25] (73.center) to (74.center);
		\draw [in=0, out=180, looseness=1.25] (75.center) to (74.center);
		\draw [in=-180, out=45, looseness=1.25] (78.center) to (76.center);
		\draw [in=-45, out=180, looseness=1.25] (77.center) to (78.center);
		\draw [in=-180, out=0, looseness=1.25] (79.center) to (78.center);
		\draw [in=0, out=135, looseness=1.25] (82.center) to (80.center);
		\draw [in=-135, out=0, looseness=1.25] (81.center) to (82.center);
		\draw (84.center) to (85.center);
		\draw (86.center) to (87.center);
		\draw (88.center) to (89.center);
		\draw (90.center) to (91.center);
		\draw (93.center) to (94.center);
		\draw (95.center) to (96.center);
		\draw (100.center) to (101.center);
		\draw [style=pointille] (102.center) to (103.center);
		\draw [style=pointille] (104.center) to (105.center);
		\draw [style=pointille] (103.center) to (105.center);
		\draw [style=pointille] (102.center) to (104.center);
		\draw [in=-180, out=45, looseness=1.25] (108.center) to (106.center);
		\draw [in=-45, out=180, looseness=1.25] (107.center) to (108.center);
		\draw [in=-180, out=0, looseness=1.25] (109.center) to (108.center);
		\draw (114.center) to (115.center);
		\draw [in=-180, out=45, looseness=1.25] (118.center) to (116.center);
		\draw [in=-45, out=180, looseness=1.25] (117.center) to (118.center);
		\draw [in=-180, out=0, looseness=1.25] (119.center) to (118.center);
		\draw (121.center) to (122.center);
		\draw [in=0, out=135, looseness=1.25] (125.center) to (123.center);
		\draw [in=-135, out=0, looseness=1.25] (124.center) to (125.center);
		\draw [in=0, out=180, looseness=1.25] (126.center) to (125.center);
		\draw (128.center) to (129.center);
		\draw [in=0, out=135, looseness=1.25] (132.center) to (130.center);
		\draw [in=-135, out=0, looseness=1.25] (131.center) to (132.center);
		\draw [in=0, out=180, looseness=1.25] (133.center) to (132.center);
		\draw (135.center) to (136.center);
	\end{pgfonlayer}
\end{tikzpicture}

%% file: figures/arrowbox.tikz
\begin{tikzpicture}
	\begin{pgfonlayer}{nodelayer}
		\node [style=none] (0) at (1, 1) {};
		\node [style=none] (1) at (1, -1) {};
		\node [style=none] (2) at (3, 1) {};
		\node [style=none] (3) at (3, -1) {};
		\node [style=none] (4) at (0, 0.5) {};
		\node [style=none] (5) at (0, -0.5) {};
		\node [style=none] (6) at (3, 0.5) {};
		\node [style=none] (7) at (4, 0.5) {};
		\node [style=none] (8) at (1, -0.5) {};
		\node [style=none] (9) at (1, 0.5) {};
		\node [style=none] (10) at (2, 0) {$[f]$};
		\node [style=none] (17) at (4, -0.5) {};
		\node [style=none] (18) at (3, -0.5) {};
		\node [style=none] (24) at (0.5, 0) {$\bvdots$};
		\node [style=none] (25) at (3.5, 0) {$\bvdots$};
	\end{pgfonlayer}
	\begin{pgfonlayer}{edgelayer}
		\draw (0.center) to (1.center);
		\draw (1.center) to (3.center);
		\draw (3.center) to (2.center);
		\draw (2.center) to (0.center);
		\draw (4.center) to (9.center);
		\draw (5.center) to (8.center);
		\draw (6.center) to (7.center);
		\draw (17.center) to (18.center);
	\end{pgfonlayer}
\end{tikzpicture}

%% file: figures/arrowboxex.tikz
\begin{tikzpicture}
	\begin{pgfonlayer}{nodelayer}
		\node [style=none] (0) at (1, 1.5) {};
		\node [style=none] (1) at (1, -0.5) {};
		\node [style=none] (2) at (3, 1.5) {};
		\node [style=none] (3) at (3, -0.5) {};
		\node [style=none] (4) at (0, 1) {};
		\node [style=none] (5) at (0, 0) {};
		\node [style=none] (6) at (3, 1) {};
		\node [style=none] (7) at (4, 1) {};
		\node [style=none] (8) at (1, 0) {};
		\node [style=none] (9) at (1, 1) {};
		\node [style=none] (10) at (2, 0.5) {$[f]$};
		\node [style=none] (17) at (4, 0) {};
		\node [style=none] (18) at (3, 0) {};
		\node [style=none] (19) at (4, -1) {};
		\node [style=none] (20) at (4, -3) {};
		\node [style=none] (21) at (6, -1) {};
		\node [style=none] (22) at (6, -3) {};
		\node [style=none] (23) at (3, -1.5) {};
		\node [style=none] (24) at (3, -2.5) {};
		\node [style=none] (25) at (6, -1.5) {};
		\node [style=none] (26) at (7, -1.5) {};
		\node [style=none] (27) at (4, -2.5) {};
		\node [style=none] (28) at (4, -1.5) {};
		\node [style=none] (29) at (5, -2) {$[g]$};
		\node [style=none] (30) at (7, -2.5) {};
		\node [style=none] (31) at (6, -2.5) {};
		\node [style=none] (32) at (11, 1) {};
		\node [style=none] (33) at (11, -1) {};
		\node [style=none] (34) at (13, 1) {};
		\node [style=none] (35) at (13, -1) {};
		\node [style=none] (36) at (10, 0.5) {};
		\node [style=none] (37) at (10, -0.5) {};
		\node [style=none] (38) at (13, 0.5) {};
		\node [style=none] (39) at (14, 0.5) {};
		\node [style=none] (40) at (11, -0.5) {};
		\node [style=none] (41) at (11, 0.5) {};
		\node [style=none] (42) at (12, 0) {$[h]$};
		\node [style=none] (43) at (14, -0.5) {};
		\node [style=none] (44) at (13, -0.5) {};
		\node [style=new style 0] (47) at (9, 0.5) {};
		\node [style=new style 0] (48) at (9, -0.5) {};
		\node [style=new style 0] (49) at (2, -2) {};
		\node [style=none] (50) at (0, -1) {};
		\node [style=none] (51) at (15, 0) {$=$};
		\node [style=none] (53) at (17, 1.5) {};
		\node [style=none] (54) at (17, -1.5) {};
		\node [style=none] (55) at (21, 1.5) {};
		\node [style=none] (56) at (21, -1.5) {};
		\node [style=none] (57) at (16, 1) {};
		\node [style=none] (58) at (16, 0) {};
		\node [style=none] (59) at (21, 0.5) {};
		\node [style=none] (60) at (22, 0.5) {};
		\node [style=none] (61) at (17, 0) {};
		\node [style=none] (62) at (17, 1) {};
		\node [style=none] (63) at (19, 0) {$[h\circ (f\otimes g)]$};
		\node [style=none] (64) at (22, -0.5) {};
		\node [style=none] (65) at (21, -0.5) {};
		\node [style=none] (66) at (16, -1) {};
		\node [style=none] (67) at (17, -1) {};
	\end{pgfonlayer}
	\begin{pgfonlayer}{edgelayer}
		\draw (0.center) to (1.center);
		\draw (1.center) to (3.center);
		\draw (3.center) to (2.center);
		\draw (2.center) to (0.center);
		\draw (4.center) to (9.center);
		\draw (5.center) to (8.center);
		\draw (6.center) to (7.center);
		\draw (17.center) to (18.center);
		\draw (19.center) to (20.center);
		\draw (20.center) to (22.center);
		\draw (22.center) to (21.center);
		\draw (21.center) to (19.center);
		\draw (23.center) to (28.center);
		\draw (24.center) to (27.center);
		\draw (25.center) to (26.center);
		\draw (30.center) to (31.center);
		\draw (32.center) to (33.center);
		\draw (33.center) to (35.center);
		\draw (35.center) to (34.center);
		\draw (34.center) to (32.center);
		\draw (36.center) to (41.center);
		\draw (37.center) to (40.center);
		\draw (38.center) to (39.center);
		\draw (43.center) to (44.center);
		\draw [in=180, out=0, looseness=1.25] (7.center) to (47);
		\draw [in=-135, out=0] (26.center) to (47);
		\draw (47) to (36.center);
		\draw [in=-120, out=0] (30.center) to (48);
		\draw [in=0, out=180, looseness=1.25] (48) to (17.center);
		\draw (48) to (37.center);
		\draw [in=-180, out=45, looseness=1.25] (49) to (23.center);
		\draw [in=-45, out=180, looseness=1.25] (24.center) to (49);
		\draw [in=-180, out=0, looseness=1.25] (50.center) to (49);
		\draw (53.center) to (54.center);
		\draw (54.center) to (56.center);
		\draw (56.center) to (55.center);
		\draw (55.center) to (53.center);
		\draw (57.center) to (62.center);
		\draw (58.center) to (61.center);
		\draw (59.center) to (60.center);
		\draw (64.center) to (65.center);
		\draw (66.center) to (67.center);
	\end{pgfonlayer}
\end{tikzpicture}

%% file: figures/ambadj.tikz
\begin{tikzpicture}
	\begin{pgfonlayer}{nodelayer}
		\node [style=none] (52) at (5, 0) {$\textbf{Prop}$};
		\node [style=none] (53) at (0, 0) {$\textbf{SSMC}$};
		\node [style=none] (54) at (1, 0.5) {};
		\node [style=none] (55) at (4, 0.5) {};
		\node [style=none] (56) at (4, -0.5) {};
		\node [style=none] (57) at (1, -0.5) {};
		\node [style=none] (58) at (2.5, 1.5) {$\textbf{Prop}$};
		\node [style=none] (59) at (2.5, -1.5) {$\overline{\phantom{ }\cdot\phantom{ }}$};
		\node [style=none] (60) at (2.5, 0) {\rotatebox{90}{$\vdash$}};
	\end{pgfonlayer}
	\begin{pgfonlayer}{edgelayer}
		\draw [style=arrow, bend left] (54.center) to (55.center);
		\draw [style=arrow, bend left] (56.center) to (57.center);
	\end{pgfonlayer}
\end{tikzpicture}

%% file: figures/unitcounitdef.tikz
\begin{tikzpicture}
	\begin{pgfonlayer}{nodelayer}
		\node [style=none] (207) at (18, -1.5) {$\mathbf{Prop}$};
		\node [style=none] (222) at (19, 3.5) {$\quad\langle\_\rangle: S \Rightarrow id_{\mathbf{Prop}}\quad$};
		\node [style=none] (224) at (7, 3.5) {$\quad[\_]: id_{\mathbf{SSMC}} \Rightarrow \overline{\mathbf{Prop}(\_)}$};
		\node [style=none] (229) at (6, 0.5) {};
		\node [style=none] (230) at (8, 0.5) {};
		\node [style=none] (231) at (7, -0.5) {};
		\node [style=none] (234) at (19, 1.5) {$id_{\mathbf{Prop}}$};
		\node [style=none] (247) at (8, 1.5) {$\mathbf{Prop}$};
		\node [style=none] (249) at (7, -1.5) {$id_{\mathbf{SSMC}}$};
		\node [style=none] (251) at (20, -1.5) {$\overline{\phantom{ }\cdot\phantom{ }}$};
		\node [style=none] (252) at (6, 1.5) {$\overline{\phantom{ }\cdot\phantom{ }}$};
		\node [style=none] (253) at (18, -0.5) {};
		\node [style=none] (254) at (20, -0.5) {};
		\node [style=none] (255) at (19, 0.5) {};
	\end{pgfonlayer}
	\begin{pgfonlayer}{edgelayer}
		\draw [bend right=45] (229.center) to (231.center);
		\draw [bend right=45] (231.center) to (230.center);
		\draw [bend left=45] (253.center) to (255.center);
		\draw [bend left=45] (255.center) to (254.center);
	\end{pgfonlayer}
\end{tikzpicture}

%% file: figures/unitcounitequations.tikz
\begin{tikzpicture}
	\begin{pgfonlayer}{nodelayer}
		\node [style=none] (59) at (0, 0) {};
		\node [style=none] (60) at (1, 0) {};
		\node [style=none] (61) at (0.5, 0.5) {};
		\node [style=none] (62) at (1, 0) {};
		\node [style=none] (63) at (2, 0) {};
		\node [style=none] (64) at (1.5, -0.5) {};
		\node [style=none] (71) at (2.75, 0) {$=$};
		\node [style=none] (72) at (2, 1) {};
		\node [style=none] (73) at (0, -1) {};
		\node [style=none] (74) at (3.5, 1) {};
		\node [style=none] (75) at (3.5, -1) {};
		\node [style=none] (85) at (9.5, 1) {};
		\node [style=none] (86) at (9.5, -1) {};
		\node [style=none] (87) at (10.25, 0) {$=$};
		\node [style=none] (88) at (13, 0) {};
		\node [style=none] (89) at (12, 0) {};
		\node [style=none] (90) at (12.5, 0.5) {};
		\node [style=none] (91) at (12, 0) {};
		\node [style=none] (92) at (11, 0) {};
		\node [style=none] (93) at (11.5, -0.5) {};
		\node [style=none] (94) at (11, 1) {};
		\node [style=none] (95) at (13, -1) {};
	\end{pgfonlayer}
	\begin{pgfonlayer}{edgelayer}
		\draw [bend left=45] (59.center) to (61.center);
		\draw [bend left=45] (61.center) to (60.center);
		\draw [bend right=45] (62.center) to (64.center);
		\draw [bend right=45] (64.center) to (63.center);
		\draw (73.center) to (59.center);
		\draw (72.center) to (63.center);
		\draw (74.center) to (75.center);
		\draw (85.center) to (86.center);
		\draw [bend right=45] (88.center) to (90.center);
		\draw [bend right=45] (90.center) to (89.center);
		\draw [bend left=45] (91.center) to (93.center);
		\draw [bend left=45] (93.center) to (92.center);
		\draw (95.center) to (88.center);
		\draw (94.center) to (92.center);
	\end{pgfonlayer}
\end{tikzpicture}

%% file: figures/monadcomonaddef.tikz
\begin{tikzpicture}
	\begin{pgfonlayer}{nodelayer}
		\node [style=none] (59) at (4.5, -0.5) {};
		\node [style=none] (60) at (5.5, -0.5) {};
		\node [style=none] (61) at (5, 0) {};
		\node [style=none] (62) at (5.5, -0.5) {};
		\node [style=none] (63) at (6.5, -1.5) {};
		\node [style=none] (71) at (2, 0) {$\doteq$};
		\node [style=none] (72) at (6.5, -0.5) {};
		\node [style=none] (73) at (4.5, -1.5) {};
		\node [style=white] (146) at (0.5, 0) {};
		\node [style=none] (147) at (0.5, 1) {};
		\node [style=none] (148) at (0, -1) {};
		\node [style=none] (149) at (1, -1) {};
		\node [style=white] (154) at (13, -0.5) {};
		\node [style=none] (155) at (13, 0.5) {};
		\node [style=none] (159) at (14, 0) {$\doteq$};
		\node [style=none] (161) at (15, -0.25) {};
		\node [style=none] (162) at (15, -0.25) {};
		\node [style=none] (163) at (16, -0.25) {};
		\node [style=none] (164) at (15.5, -0.75) {};
		\node [style=none] (165) at (16, 0.75) {};
		\node [style=none] (166) at (15, 0.75) {};
		\node [style=none] (174) at (5.5, -1.5) {};
		\node [style=none] (175) at (3.5, -0.5) {};
		\node [style=none] (176) at (4.5, 2) {};
		\node [style=none] (178) at (4.5, 2) {};
		\node [style=none] (179) at (3.5, -1.5) {};
		\node [style=none] (180) at (4.5, 1) {};
		\node [style=none] (181) at (5.5, 2) {};
		\node [style=none] (182) at (5.5, 2) {};
		\node [style=none] (183) at (5.5, 1) {};
		\node [style=none] (205) at (0, -1.75) {$S$};
		\node [style=none] (206) at (1, -1.75) {$S$};
		\node [style=none] (207) at (0.5, 1.75) {$S$};
		\node [style=none] (210) at (13, 1.25) {$S$};
	\end{pgfonlayer}
	\begin{pgfonlayer}{edgelayer}
		\draw [bend left=45] (59.center) to (61.center);
		\draw [bend left=45] (61.center) to (60.center);
		\draw (73.center) to (59.center);
		\draw (72.center) to (63.center);
		\draw [in=-135, out=90, looseness=1.25] (148.center) to (146);
		\draw [in=90, out=-45, looseness=1.25] (146) to (149.center);
		\draw (146) to (147.center);
		\draw (154) to (155.center);
		\draw [bend right=45] (162.center) to (164.center);
		\draw [bend right=45] (164.center) to (163.center);
		\draw (165.center) to (163.center);
		\draw (166.center) to (162.center);
		\draw (62.center) to (174.center);
		\draw (179.center) to (175.center);
		\draw (178.center) to (180.center);
		\draw (182.center) to (183.center);
		\draw [in=90, out=-90, looseness=1.25] (183.center) to (72.center);
		\draw [in=90, out=-90, looseness=1.25] (180.center) to (175.center);
	\end{pgfonlayer}
\end{tikzpicture}

%% file: figures/adjstrict.tikz
\begin{tikzpicture}
	\begin{pgfonlayer}{nodelayer}
		\node [style=none] (52) at (5, 0) {$\textbf{SSMC}$};
		\node [style=none] (53) at (0, 0) {$\textbf{SMC}$};
		\node [style=none] (54) at (1, 0.5) {};
		\node [style=none] (55) at (4, 0.5) {};
		\node [style=none] (56) at (4, -0.5) {};
		\node [style=none] (57) at (1, -0.5) {};
		\node [style=none] (58) at (2.5, 1.5) {$\textbf{Strict}$};
		\node [style=none] (59) at (2.5, -1.5) {$\overline{\phantom{ }\cdot\phantom{ }}$};
		\node [style=none] (60) at (2.5, 0) {\rotatebox{90}{$\vdash$}};
	\end{pgfonlayer}
	\begin{pgfonlayer}{edgelayer}
		\draw [style=arrow, bend left] (54.center) to (55.center);
		\draw [style=arrow, bend left] (56.center) to (57.center);
	\end{pgfonlayer}
\end{tikzpicture}

%% file: figures/adjprop.tikz
\begin{tikzpicture}
	\begin{pgfonlayer}{nodelayer}
		\node [style=none] (52) at (5, 0) {$\textbf{Prop}$};
		\node [style=none] (53) at (0, 0) {$\textbf{SMC}$};
		\node [style=none] (54) at (1, 0.5) {};
		\node [style=none] (55) at (4, 0.5) {};
		\node [style=none] (56) at (4, -0.5) {};
		\node [style=none] (57) at (1, -0.5) {};
		\node [style=none] (58) at (2.5, 1.5) {$\textbf{Prop}$};
		\node [style=none] (59) at (2.5, -1.5) {$\overline{\phantom{ }\cdot\phantom{ }}$};
		\node [style=none] (60) at (2.5, 0) {\rotatebox{90}{$\vdash$}};
		\node [style=none] (61) at (-1, 0.5) {};
		\node [style=none] (62) at (-1, -0.5) {};
		\node [style=none] (63) at (6, 0.5) {};
		\node [style=none] (64) at (6, -0.5) {};
		\node [style=none] (65) at (-3.5, 0) {$\textbf{Strict}$};
		\node [style=none] (66) at (8, 0) {$S$};
	\end{pgfonlayer}
	\begin{pgfonlayer}{edgelayer}
		\draw [style=arrow, bend left] (54.center) to (55.center);
		\draw [style=arrow, bend left] (56.center) to (57.center);
		\draw [style=arrow, bend left=135, looseness=5.50] (62.center) to (61.center);
		\draw [style=arrow, bend right=135, looseness=5.50] (64.center) to (63.center);
	\end{pgfonlayer}
\end{tikzpicture}

%% file: figures/naturalityunitscounits.tikz
\begin{tikzpicture}
	\begin{pgfonlayer}{nodelayer}
		\node [style=none] (0) at (0, 1.5) {$\overline{\mathbf{Prop}(\mathbf{C})}$};
		\node [style=none] (1) at (1, 4) {$\mathbf{C}$};
		\node [style=none] (4) at (2, 4) {};
		\node [style=none] (5) at (5, 4) {};
		\node [style=none] (6) at (6, 4) {$\mathbf{D}$};
		\node [style=none] (7) at (1.5, 1.5) {};
		\node [style=none] (8) at (5.5, 1.5) {};
		\node [style=none] (9) at (1, -3.5) {};
		\node [style=none] (10) at (6, -3.5) {};
		\node [style=none] (11) at (0, -2.5) {};
		\node [style=none] (12) at (0, 0.5) {};
		\node [style=none] (13) at (7, -2.5) {};
		\node [style=none] (14) at (7, 0.5) {};
		\node [style=none] (19) at (0, 4) {$F:$};
		\node [style=none] (21) at (7, 1.5) {$\overline{\mathbf{Prop}(\mathbf{D})}$};
		\node [style=none] (22) at (0, -3.5) {$\mathbf{C}$};
		\node [style=none] (23) at (7, -3.5) {$\mathbf{D}$};
		\node [style=none] (24) at (14, 1.5) {$S(\mathbf{P})$};
		\node [style=none] (25) at (15, 4) {$\mathbf{P}$};
		\node [style=none] (26) at (16, 4) {};
		\node [style=none] (27) at (19, 4) {};
		\node [style=none] (28) at (20, 4) {$\mathbf{Q}$};
		\node [style=none] (29) at (15, 1.5) {};
		\node [style=none] (30) at (20, 1.5) {};
		\node [style=none] (31) at (15, -3.5) {};
		\node [style=none] (32) at (20, -3.5) {};
		\node [style=none] (37) at (14, 0.5) {};
		\node [style=none] (38) at (14, -2.5) {};
		\node [style=none] (39) at (21, 0.5) {};
		\node [style=none] (40) at (21, -2.5) {};
		\node [style=none] (41) at (14, 4) {$G:$};
		\node [style=none] (43) at (14, -3.5) {$\mathbf{P}$};
		\node [style=none] (44) at (21, -3.5) {$\mathbf{Q}$};
		\node [style=none] (45) at (21, 1.5) {$S(\mathbf{Q})$};
		\node [style=none] (47) at (17.5, 2.25) {$S(G)$};
		\node [style=none] (48) at (3.5, 2.25) {$\overline{\mathbf{Prop}(F)}$};
		\node [style=none] (49) at (17.5, -2.75) {$G$};
		\node [style=none] (50) at (3.5, -2.75) {$F$};
		\node [style=none] (51) at (-1, -1) {$[\_]_{\mathbf{C}}$};
		\node [style=none] (52) at (8, -1) {$[\_]_{\mathbf{D}}$};
		\node [style=none] (57) at (13, -1) {$\langle\_\rangle_{\mathbf{P}}$};
		\node [style=none] (58) at (22, -1) {$\langle\_\rangle_{\mathbf{Q}}$};
	\end{pgfonlayer}
	\begin{pgfonlayer}{edgelayer}
		\draw [style=arrow] (4.center) to (5.center);
		\draw [style=arrow] (7.center) to (8.center);
		\draw [style=arrow] (9.center) to (10.center);
		\draw [style=arrow] (11.center) to (12.center);
		\draw [style=arrow] (13.center) to (14.center);
		\draw [style=arrow] (26.center) to (27.center);
		\draw [style=arrow] (29.center) to (30.center);
		\draw [style=arrow] (31.center) to (32.center);
		\draw [style=arrow] (37.center) to (38.center);
		\draw [style=arrow] (39.center) to (40.center);
	\end{pgfonlayer}
\end{tikzpicture}